\documentclass{siamart1116}

\usepackage[margin=1in]{geometry}
\usepackage{graphicx}
\usepackage{color}
\usepackage{amssymb,amsmath,mathtools}
\usepackage{stmaryrd}

\usepackage{bm}
\usepackage{epstopdf}
\usepackage{amsfonts}
\usepackage{paralist}
\usepackage{algorithm}
\usepackage{algcompatible}
\usepackage{pstool}
\usepackage{wrapfig}
\usepackage{rotating}

\usepackage{tikz}
\usepackage{pgfplots}
\usepackage{pgfplotstable, filecontents}
\pgfplotsset{compat=1.9}
\usetikzlibrary{pgfplots.groupplots}
\usepgfplotslibrary{fillbetween}
\usetikzlibrary{calc,fit,matrix,arrows,automata,positioning,shapes}
\usetikzlibrary{arrows.meta}
\usepackage{pgfplotstable, booktabs}
\pgfplotsset{select coords between index/.style 2 args={
    x filter/.code={
        \ifnum\coordindex<#1\fi
        \ifnum\coordindex>#2\fi
    }
}}

\ifpdf
  \DeclareGraphicsExtensions{.eps,.pdf,.png,.jpg}
\else
  \DeclareGraphicsExtensions{.eps}
\fi

\numberwithin{theorem}{section}

\newcommand{\TheTitle}{An efficient, globally convergent method for
                       optimization under uncertainty using adaptive
                       model reduction and sparse grids}
\newcommand{\TheShortTitle}{Adaptive model reduction, stochastic collocation}
\newcommand{\TheAuthors}{M.\ J.\ Zahr, K.\ T.\ Carlberg, D.\ P.\ Kouri}
\headers{\TheShortTitle}{\TheAuthors}

\title{{\TheTitle}\thanks{Submitted to the editors DATE.}}

\author{
  Matthew J. Zahr%
  \thanks{Department of Mathematics, Lawrence Berkeley National Laboratory,
          University of California, Berkeley, Berkeley, CA
          (\email{mjzahr@lbl.gov},
           \url{http://math.lbl.gov/\textasciitilde mjzahr});
          Department of Aerospace and Mechanical Engineering, University of
          Notre Dame, Notre Dame, IN (\email{mzahr@nd.edu})}
  \and
  Kevin T.\ Carlberg%
  \thanks{Sandia National Laboratories, Livermore, CA
          (\email{ktcarlb@sandia.gov})}
  \and
  Drew P.\ Kouri%
  \thanks{Sandia National Laboratories, Albuquerque, NM
          (\email{dpkouri@sandia.gov}) \newline
     Sandia National Laboratories is a multimission laboratory
     managed and operated by National Technology and Engineering
     Solutions of Sandia, LLC., a wholly owned subsidiary of
     Honeywell International, Inc., for the U.S.\ Department of
     Energy’s National Nuclear Security Administration under
     contract DE-NA0003525.
     This paper describes objective technical results and analysis. Any
     subjective views or opinions that might be expressed in the paper
     do not necessarily represent the views of the U.S.\ Department of
     Energy or the United States Government.}
}

\usepackage{amsopn}

\newtheorem{assume}{Assumption}

\newtheorem{prop}{Proposition}

\newtheorem*{rem}{Remark}

\newcommand{\func}[3]{\ensuremath{#1 : #2 \rightarrow #3}}

\newcommand{\norm}[1]{\ensuremath{\left\| #1 \right\|}}
\newcommand{\optunc}[2]{\underset{#1}{\mathrm{minimize}} ~~ #2}
\newcommand{\optcona}[3]{
\begin{aligned}
& \underset{#1}{\text{minimize}}
& & #2 \\
& \text{subject to} & & #3
\end{aligned}}

\newcommand{\pder}[2]{\ensuremath{\frac{\partial #1}{\partial #2}}} %1st partial derivative
\newcommand{\Bcal}{\ensuremath{\mathcal{B}}}

\newcommand{\Ecal}{\ensuremath{\mathcal{E}}}

\newcommand{\Ical}{\ensuremath{\mathcal{I}}}
\newcommand{\Jcal}{\ensuremath{\mathcal{J}}}

\newcommand{\Ncal}{\ensuremath{\mathcal{N}}}

\newcommand{\Rcal}{\ensuremath{\mathcal{R}}}

\newcommand{\Ucal}{\ensuremath{\mathcal{U}}}

\newcommand{\Zcal}{\ensuremath{\mathcal{Z}}}

\newcommand{\Ebb}{\ensuremath{\mathbb{E}}}

\newcommand{\Nbb}{\ensuremath{\mathbb{N} }}
\newcommand{\Qbb}{\ensuremath{\mathbb{Q} }}
\newcommand{\Rbb}{\ensuremath{\mathbb{R} }}
\newcommand{\Rbbstar}{\ensuremath{\mathbb{R}_\star}}

\newcommand\Dbm{{\ensuremath{\bm{D}}}}

\newcommand\Ibm{{\ensuremath{\bm{I}}}}

\newcommand\ebm{{\ensuremath{\bm{e}}}}

\newcommand\gbm{{\ensuremath{\bm{g}}}}

\newcommand\ibm{{\ensuremath{\bm{i}}}}

\newcommand\kbm{{\ensuremath{\bm{k}}}}

\newcommand\qbm{{\ensuremath{\bm{q}}}}
\newcommand\rbm{{\ensuremath{\bm{r}}}}

\newcommand\ubm{{\ensuremath{\bm{u}}}}

\newcommand\xbm{{\ensuremath{\bm{x}}}}
\newcommand\ybm{{\ensuremath{\bm{y}}}}
\newcommand\zbm{{\ensuremath{\bm{z}}}}

\newcommand\lambdabold{{\ensuremath{\boldsymbol{\lambda}}}}

\newcommand\etabold{{\ensuremath{\boldsymbol{\eta}}}}

\newcommand\mubold{{\ensuremath{\boldsymbol{\mu}}}}

\newcommand\Phibold{{\ensuremath{\boldsymbol{\Phi}}}}

\newcommand\Thetabold{{\ensuremath{\boldsymbol{\Theta}}}}

\newcommand\Psibold{{\ensuremath{\boldsymbol{\Psi}}}}
\newcommand\Xibold{{\ensuremath{\boldsymbol{\Xi}}}}

\newcommand\zerobold{\ensuremath{\mathbf{0}}}

\newcommand\Range{\ensuremath{\mathrm{Ran}}}

\usepackage{etoolbox}
\newbool{fastcompile}
\setbool{fastcompile}{false}

\usepackage{xifthen}
\newcommand{\nU}{{n_{\ubm}}}
\newcommand{\kU}{{k_{\ubm}}}
\newcommand{\ny}{{n_{\ybm}}}
\newcommand{\nmu}{{n_{\mubold}}}
\newcommand{\nz}{n_{\zbm}}
\newcommand{\nzVerbose}{\ny+\nmu}
\newcommand{\dpdesol}{\ubm}

\newcommand{\dpdesolSol}{\ubm_\star}

\newcommand{\dadjvar}{\lambdabold}
\newcommand{\dadjvarSol}{\lambdabold_\star}
\newcommand{\dpderes}{\rbm}
\newcommand{\dadjres}{\rbm^\lambdabold}
\newcommand{\dadjgrad}{\gbm^\lambdabold}
\newcommand{\dstochvar}{\ybm}

\newcommand{\dparamvar}{\mubold}

\newcommand{\dcombvar}{\zbm}
\newcommand{\dstochsp}[1][]{\ifthenelse{\isempty{#1}}{\Xibold}{\Xi_{#1}}}
\newcommand{\probdens}[1][]{\ifthenelse{\isempty{#1}}{\rho}{\rho_{#1}}}

\newcommand{\lrob}{\Psibold}
\newcommand{\rrob}{\Phibold}
\newcommand{\dromsol}{\qbm}
\newcommand{\dromsolSol}{\qbm_\star}
\newcommand{\dromadj}{\etabold}
\newcommand{\dromadjSol}{\etabold_\star}
\newcommand{\dromadjmr}{\hat\dromadj}
\newcommand{\dromadjmrSol}{\hat\dromadj_\star}

\newcommand{\obj}{J}
\newcommand{\qoi}{f}
\newcommand{\qoirstr}{F}
\newcommand{\rqoirstr}{F_r}
\newcommand{\riskm}{\Rcal}

\newcommand{\trrad}{\Delta}
\newcommand{\aprxmod}{m}
\newcommand{\aprxmodobj}{\psi}
\newcommand{\objbnd}{\theta}
\newcommand{\gradbnd}{\varphi}

\newcommand{\eqnref}[1]{Eq.~\eqref{#1}}

\ifpdf
\hypersetup{
  pdftitle={\TheTitle},
  pdfauthor={\TheAuthors}
}
\fi

\begin{document}

\maketitle

% REQUIRED: single paragraph, only inline equations, no citations to paper's references
\begin{abstract}
This work introduces a new method to efficiently solve optimization problems
constrained by partial differential equations (PDEs) with uncertain
coefficients. The method leverages two sources of inexactness that trade
accuracy for speed: (1) stochastic collocation based on dimension-adaptive
sparse grids (SGs), which approximates the stochastic objective function
with a limited number of quadrature nodes, and (2) projection-based
reduced-order models (ROMs), which generate efficient approximations to
PDE solutions.  These two sources of inexactness lead to inexact objective function
and gradient evaluations, which are managed by a trust-region method that
guarantees global convergence by adaptively refining the sparse grid
and reduced-order model until a proposed error indicator drops below
a tolerance specified by trust-region convergence theory. A key feature
of the proposed method is that the error indicator---which accounts for
errors incurred by both
the sparse grid and reduced-order model---must be only an
\emph{asymptotic} error bound, i.e., a bound that holds up to an arbitrary
constant that need not be computed.  This enables the method to be applicable
to a wide range of problems, including those where sharp, computable error
bounds are not available; this distinguishes the proposed method from
previous works. Numerical experiments performed on a model
problem from optimal flow control under uncertainty verify global
convergence of the method and demonstrate the method's ability to
outperform previously proposed alternatives.
\end{abstract}

% REQUIRED

% REQUIRED

% Usage: \cref{X} for "section #"/"theorem #"/etc (\Cref for capital), \ref{X} for "#"
\section{Introduction}

Optimization problems constrained by parametrized systems of nonlinear
equations---which usually result from the discretization of a system of partial
differential equations---arise in nearly every branch of engineering and
science and in many different contexts, including design, control, and data
assimilation. To compute solutions that are robust with respect to 
uncertainties that affect physical systems, these uncertainties must be
incorporated within the optimization formulation. This leads to a coupling
between optimization and uncertainty quantification, and can introduce a
significant computational burden when the nonlinear system and stochastic
space are high-dimensional.

Algorithms proposed in Refs.~\cite{kouri2013trust, kouri2014inexact} have
shown promise in reducing the cost of optimization under uncertainty
compared to previously proposed methods. These algorithms enable computational
efficiency via adaptive sparse-grid stochastic collocation and a
practical trust-region framework to manage the resulting models. While these
algorithms constitute important steps toward enabling large-scale optimization
under uncertainty, they often require thousands of solutions to the
parametrized nonlinear system, which is prohibitively expensive for
large-scale systems. Indeed, the primary cost associated with the method
proposed in Ref.~\cite{kouri2014inexact} arises from the need to compute
primal and adjoint solutions associated with the parametrized nonlinear system
at the prescribed stochastic collocation nodes; these solutions are needed to
define the trust-region model and associated error indicators. The goal of this
paper is to mitigate this computational burden, which we achieve by allowing
for inexact solutions of the parametrized nonlinear system at collocation
nodes. In particular, we use minimum-residual, projection-based reduced-order
models to efficiently approximate solutions to the parametrized nonlinear
system, i.e., high-dimensional model (HDM), and bootstrap the first-order
convergence theory introduced in Ref.~\cite{kouri2014inexact} to ensure
global convergence of the resulting method.
%\DPK{The algorithms referenced above do not technically require exact forward
%and adjoint solves.  For example, one could also add in adaptive FEM.  We
%chose to focus on adaptive sparse grids to simplify presentation.  We might
%want to mention that the trust-region algorithm does not require exact 
%solves, but the stated adaptive sparse grid algorithms (include algorithm
%numbers) do.}

The proposed method begins by approximating the objective function in the
optimization under uncertainty problem using sparse grids and reduced-order
models. Sparse grids provide efficient quadrature rules with natural
dimension-adaptive refinement and error estimation to approximate the integrals
that arise in the computation of risk measures \cite{artzner1999coherent};
reduced-order models cheaply approximate the solution to the parametrized
nonlinear system at the quadrature nodes.
%\DPK{Sparse grids require regularity (smoothness) from the integrand, but
%coherent risk measures are never differentiable (positive homogeneity
%kills the possibility of differentiability).}
The method solves the resulting optimization problem in a globally convergent
manner using the trust-region method proposed in Ref.~\cite{kouri2014inexact},
which allows for inexact objective and gradient evaluation, even at
trust-region centers. Global first-order convergence of this trust-region
method relies on (1) inexpensive, computable error indicators for the
objective and gradient approximations and (2) models that can be refined
until these error indicators drop below required tolerances.  The only
requirement on the error indicators is that they bound the corresponding
error \emph{up to an arbitrary constant} that need not be computed or
estimated; this provides substantial flexibility to develop
\textit{a posteriori} residual-based error indicators that are applicable
to a wide range of problems.
%\DPK{This paragraph addresses my concern with the previous paragraph.  Perhaps
%we should combine them somehow.}
To this end, we derive error indicators for the
objective function and its gradient. The gradient error indicator comprises a
combination of both primal and adjoint residual-based error
indicators for the ROM approximation \cite{zahr2015progressive, zahr2016phd}
and an approximation of the quadrature truncation error on the forward
neighbors of the sparse grid \cite{gerstner2003dimension}.

With the error indicators defined, we introduce two algorithms to construct an
anisotropic sparse grid and reduced basis, which together form the trust
region approximation models. The first algorithm constructs a
sparse grid and reduced basis such that an accuracy condition
on the \textit{gradient error indicator} \cite{kouri2013trust} is satisfied
at the trust-region center. The second algorithm constructs a (possibly
different) sparse grid and reduced basis such that an accuracy condition
on the \textit{objective error indicator} \cite{kouri2014inexact} is
satisfied. Both algorithms combine the dimension-adaptive approach proposed
in Ref.~\cite{gerstner2003dimension} to construct an anisotropic
sparse grid  with a variant of the greedy method proposed in
Refs.~\cite{patera2007reduced,rozza2008reduced} to construct a reduced
basis into a single nested algorithm that constructs both simultaneously.
Here, the outer loop
refines the sparse grid by adding the index set from the set of forward
neighbors that contributes most to the quadrature truncation error for the ROM
gradient.  A direct search for the index set that maximizes this error is
inexpensive because it requires only ROM evaluations. Once the sparse grid is
refined, the inner loop greedily constructs the reduced basis by
sampling the HDM at sparse grid collocation points where
the primal and adjoint residual-based error indicators are maximized. To ensure
the reduced basis can accurately approximate the primal and adjoint solutions,
it is constructed from both primal and dual snapshots. A crucial ingredient in
this setting is the use of \emph{minimum-residual} primal and adjoint
reduced-order models, as this approach
\begin{inparaenum}[(1)]
 \item ensures that the reduced-order-model solution 
       minimizes the objective and gradient error indicators
       over the associated trial subspaces,
 \item guarantees the greedy method based on the primal and adjoint
       residual-based indicators terminates with a reduced basis that
       satisfies the trust-region global convergence conditions, and
 \item ensures the presence of adjoint solutions in the primal basis does
       not degrade the primal approximation and vice versa, which is not
       guaranteed in general \cite{arian2000trust, eftang2013approximation}.
\end{inparaenum}

Chen et al.\ have related work \cite{chen2014weighted,
                                     chen2013multilevel, chen2013weighted}
that also
integrates reduced-order models and sparse grids for uncertainty quantification
and stochastic optimal control. Unlike the method proposed in this work;
however, their approach adheres to an \emph{offline--online} decomposition
that constructs a sparse grid and reduced-order model in a computationally
expensive offline stage and subsequently deploys them in a computationally
inexpensive online stage. Their methods are equipped with rigorous
\textit{a posteriori} error bounds and convergence guarantees, provided the
governing equations are linear elliptic PDEs.
These methods were extended to include adaptation of the sparse grid
and reduced-order model in \cite{chen2015new} and generalized to a class
of linear and nonlinear problems in the context of Bayesian inversion in
\cite{chen2015sparse, chen2016sparse}.
By breaking the offline--online
decomposition and globalizing the sparse-grid--reduced-order-model (SG--ROM)
approximation model with a practical trust-region method
\cite{kouri2014inexact}, we can establish global convergence for a
significantly wider range of governing equations, including nonlinear problems.

There are a number of other works that have combined reduced-order, or more
generally, surrogate models and/or sparse grids, to accelerate optimization under
uncertainty
\cite{maute2009reduced, yang2017algorithms, royset2017risk,
      ZZou_DPKouri_WAquino_2018a} and
uncertainty quantification
\cite{torlo2018stabilized, ullmann2014pod}
and others that have used reduced-order
models to estimate risk measures, including the conditional value-at-risk
\cite{heinkenschlossCVAR, ZZou_DPKouri_WAquino_2016a}.
Ref.~\cite{maute2009reduced} used reduced-order models to efficiently
approximate the solution of stochastic linear elasticity problems and
parametrized the stochastic space using polynomial chaos. While their
method does adaptively re-sample the parameter space throughout the
optimization to promote convergence, the adaptation is heuristic
and does not necessarily guarantee global convergence.
Ref.~\cite{ZZou_DPKouri_WAquino_2018a} presents a method that uses
local reduced-order models and Monte Carlo sampling inside a globally
convergent trust region method to minimize conditional value-at-risk.
The reduced-order models are local to a Voronoi cell of the stochastic
space, which limits the method to rather low-dimensional stochastic
spaces, and built using samples of the PDE solution and its gradient with
respect to the stochastic variables.
We also mention that other works have explored combining
reduced-order models and sparse grids for generating fast approximations
for parameterized systems
\cite{peherstorfer2012model,peherstorfer2013model,xiao2015non,lin2017non}.

The remainder of the paper is organized as follows.  Section~\ref{sec:prob}
formulates the optimization under uncertainty problem governed by a large-scale
system of nonlinear equations.  Section~\ref{sec:aprxmod} introduces the
primary approximation techniques that form the foundation for our method:
stochastic collocation via anisotropic sparse grid and minimum-residual,
projection-based primal and adjoint reduced-order models.
Section~\ref{sec:trammo} presents the new adaptive algorithm for
optimization under uncertainty that manages the inexactness introduced from the
sparse grid quadrature and reduced-order-model evaluations using the
trust-region method introduced in Ref.~\cite{kouri2014inexact}.
Appendix~\ref{app:resbnd} introduces and derives residual-based error bounds
used to define the error indicators required by the trust-region convergence
theory. Numerical results presented in Section~\ref{sec:num-exp} show the
proposed method dramatically reduces the number of PDE queries required to
solve two stochastic optimal flow control problems, as compared with the
benchmark method proposed in Ref.~\cite{kouri2014inexact}.

\section{Problem formulation}
\label{sec:prob}
% Purpose of this section is formulate the PDE-constrained optimization
% under uncertainty problem, along with necessary assumptions, and its
% discretization. Also introduce adjoint method, mainly to set notation
Let $(\Omega,\mathcal{F},P)$ be a probability space: $\Omega$ denotes the
set of outcomes, $\mathcal{F}\subseteq 2^\Omega$ denotes a $\sigma$-algebra of
events and $P:\mathcal{F}\to[0,1]$ denotes a probability measure.  We consider
the high-dimensional model (HDM) to be a large-scale system of nonlinear
equations parametrized by both a random vector
$\dstochvar\in\dstochsp$---which results from the mapping
$\Omega\ni\omega\mapsto\dstochvar\in\dstochsp \subseteq\Rbb^\ny$ with Lebesgue density
$\probdens:\dstochsp\to[0,+\infty)$---and deterministic parameters
$\dparamvar \in \Rbb^{\nmu}$. The corresponding problem statement is the
following:
Given a realization of the random vector $\dstochvar\in\dstochsp$ and an
instance of the parameters $\dparamvar\in \Rbb^{\nmu}$, compute the
(primal) solution $\dpdesolSol$ satisfying
\begin{equation} \label{eqn:sppde-disc}
 \dpderes(\dpdesolSol,\,\dstochvar,\,\dparamvar) = \zerobold,
\end{equation}
where $\dpderes:(\dpdesol,\,\dstochvar,\,\dparamvar)
                \mapsto\dpderes(\dpdesol,\,\dstochvar,\,\dparamvar)$
with $\dpderes:\Rbb^\nU\times\Rbb^{\ny}\times\Rbb^{\nmu}\rightarrow\Rbb^\nU$
denotes the residual.
We assume that for every $(\dstochvar,\dparamvar)$ pair, there exists a unique (primal) solution
$\dpdesolSol = \dpdesolSol(\dstochvar,\,\dparamvar)$ satisfying
\eqnref{eqn:sppde-disc}; we further assume that
the implicit map $(\dstochvar,\,\dparamvar)\mapsto \dpdesolSol$ is continuously differentiable with respect to its arguments.
We primarily consider the case where the residual $\dpderes$ arises from the
high-fidelity discretization of a parametrized,
stochastic partial differential equation (SPDE). In most practical applications,
the dimension $\nU$ is large, which causes \eqnref{eqn:sppde-disc} to be 
computationally expensive to solve.

Let $\qoi:
\Rbb^\nU\times\Rbb^{\ny}\times\Rbb^{\nmu}\rightarrow\Rbb$ denote
a scalar-valued quantity of
interest (QoI) associated with the system (e.g., 
the integral of a quantity over
a region of the domain 
in the SPDE setting) and let
%$\qoirstr(\dstochvar,\,\dparamvar)$ 
$\qoirstr:
\Rbb^{\ny}\times\Rbb^{\nmu}\rightarrow\Rbb
$ 
denote its restriction to the
manifold of solutions to \eqnref{eqn:sppde-disc} such that
\begin{equation} \label{eqn:qoi-rstr}
 \qoirstr:(\dstochvar,\,\dparamvar) \mapsto
 \qoi(\dpdesolSol(\dstochvar,\,\dparamvar),\,\dstochvar,\,\dparamvar).
 %\qoirstr:(\dstochvar,\,\dparamvar) \mapsto
 %\qoi(\dpdesol(\dstochvar,\,\dparamvar),\,\dstochvar,\,\dparamvar).
\end{equation}
such that 
$\qoirstr(\,\cdot\,,\,\dparamvar)\in L^1_{\probdens}(\dstochsp)$ for all
$\dparamvar\in\Rbb^{\nmu}$.
Using these definitions, we consider the stochastic optimization problem 
\begin{equation} \label{eqn:sppde-opt}
 \optunc{\dparamvar \in \Rbb^\nmu}{\obj(\dparamvar)},
\end{equation}
where the objective function is a relevant \emph{risk measure}
$\func{\riskm}{L^1_{\probdens}(\dstochsp)}{\Rbb}$
\cite{artzner1999coherent}
applied to the discrete quantity of interest, i.e.,
\begin{equation} \label{eqn:qoi-rstr-rm}
 \obj:\dparamvar \mapsto \riskm(\qoirstr(\,\cdot\,,\,\dparamvar)).
\end{equation}
This work considers only the risk-neutral measure, i.e., the
expectation risk measure $\Rcal \equiv \Ebb$; however, the approach
can be extended to other smooth risk measures.

Because the parameter-space dimension $\nmu$ is potentially large,
we employ the adjoint method to compute the gradient of the
quantity of interest. The problem statement associated with the HDM adjoint
problem is: 
Given a realization of the random vector $\dstochvar\in\dstochsp$, an
instance of the parameters $\dparamvar\in \Rbb^{\nmu}$, and the primal
solution
$\dpdesolSol \in \Rbb^\nU$, satisfying Eq.~\eqref{eqn:sppde-disc}, compute
the adjoint solution $\dadjvarSol \in \Rbb^\nU$ satisfying
\begin{equation}\label{eqn:adjres}
  \dadjres(\dadjvarSol,\,\dpdesolSol,\,\dstochvar,\,\dparamvar) 
    = \zerobold,
\end{equation}
where the adjoint residual is defined as
\begin{equation}\label{eqn:adjresDef}
\dadjres:(\dadjvar,\,\dpdesol,\,\dstochvar,\,\dparamvar)\mapsto\pder{\dpderes}{\dpdesol}(\dpdesol,\,\dstochvar,\,\dparamvar)^T \dadjvar -
   \pder{\qoi}{\dpdesol}(\dpdesol,\,\dstochvar,\,\dparamvar)^T.
\end{equation}
From the primal--adjoint pair $(\dpdesolSol,\,\dadjvarSol)$ satisfying
Eqs.~\eqref{eqn:sppde-disc} and \eqref{eqn:adjres}, the
gradient of the quantity of interest can be computed as
\begin{equation} \label{eqn:adjgrad}
 \nabla_\dparamvar\qoirstr(\dstochvar,\,\dparamvar) =
   \dadjgrad(\dadjvarSol,\,\dpdesolSol,\,\dstochvar,\,\dparamvar),
\end{equation}
where operator that reconstructs the gradient from the adjoint solution is
\begin{equation}\label{eqn:adjgradDef}
\dadjgrad:(\dadjvar,\,\dpdesol,\,\dstochvar,\,\dparamvar)\mapsto
   \pder{\qoi}{\dparamvar}(\dpdesol,\,\dstochvar,\,\dparamvar) -
   \dadjvar^T\pder{\rbm}{\dparamvar}(\dpdesol,\,\dstochvar,\,\dparamvar).
\end{equation}
The gradient of the objective, i.e., $\nabla\obj$, can be can be computed
directly from the gradient $\nabla_\dparamvar\qoirstr$ and the analytical
form of the risk measure. In the risk-neutral case, this corresponds to
the expectation of the gradient of $\qoirstr$, i.e.,
\begin{equation}
 \nabla\obj(\dparamvar) =
 \Ebb\left[\nabla_\dparamvar\qoirstr(\,\cdot\,,\,\dparamvar)\right].
\end{equation}

\section{Two sources of inexactness: sparse grids and model reduction}
\label{sec:aprxmod}
In this section, we introduce the two approximation tools used in this work
to accelerate the solution of the stochastic optimization problem 
\eqref{eqn:sppde-opt}: dimension-adaptive anisotropic sparse grids and
projection-based reduced-order models. Sparse grids provide a relatively
small number of quadrature nodes (and associated weights) to enable 
the efficient evaluation of the risk measure in a moderate-dimensional
stochastic space $\dstochsp$, while reduced-order models rapidly approximate
the quantity of interest at each quadrature node. Section~\ref{sec:trammo} 
combines these approximation techniques with a trust-region method
\cite{kouri2014inexact} to yield an efficient, globally convergent method for
solving stochastic optimization problems.

\subsection{Stochastic collocation based on sparse grids}
\label{sec:aprxmod:sg}
Assuming the stochastic space $\dstochsp$ corresponds to a tensor product of
one-dimensional intervals, i.e.,
$\dstochsp = \dstochsp[1] \otimes \cdots \otimes \dstochsp[\ny]$,
$\dstochsp[k] \subset \Rbb$, $k = 1,\dots,\ny$, the sparse-grid
construction begins with the definition of a one-dimensional quadrature
rule of level $i\in\mathbb{N}$ that will be used in the $k$th dimension,
which we denote by $\Ebb_k^i$. The level $i$ indicates refinement of the
one-dimensional quadrature rule such that
\begin{equation} \label{eqn:integ1}
 \Ebb_k^i[h] \rightarrow \Ebb_k[h] \coloneqq
 \int_{\dstochsp[k]} \probdens[k](y)h(y)\,dy \quad \text{as} \
 i \rightarrow \infty
\end{equation}
for $h \in C(\dstochsp[k])$, where the probability density is a product
of one-dimensional probability densities,
$\probdens = \probdens[1]\cdots\probdens[\ny]$.
Let $\dstochsp[k]^i \subset [-1,\,1]$ be the quadrature nodes associated with
the quadrature rule $\Ebb_k^i$. From the one-dimensional quadrature rules,
the corresponding difference operators are defined as
\begin{equation} \label{eqn:diffop0}
 \Delta_k^1 \coloneqq \Ebb_k^1 \qquad \text{and} \qquad
 \Delta_k^i \coloneqq \Ebb_k^i - \Ebb_k^{i-1} \quad \text{for} \quad
 i \geq 2.
\end{equation}
The requirement in (\ref{eqn:integ1}) on the quadrature rules implies
$\Delta_k^i[h] \rightarrow 0$ as $i \rightarrow \infty$.

A multi-dimensional difference operator is constructed from a tensor product
of one-dimensional difference operators, each possibly at a different level
of refinement, i.e.,
\begin{equation}
 \Delta^\ibm \coloneqq
 \Delta_1^{i_1} \otimes \cdots \otimes
 \Delta_{\ny}^{i_{\ny}}.
\end{equation}
A multi-index $\ibm \equiv (i_1,\,\dots,\,i_{\ny})\in \Nbb^{\ny}$ is used to track the refinement level
of each one-dimensional difference operator, i.e., $i_k$ is the refinement
level of the difference operator in dimension $k$. From the multi-dimensional
difference operator, a quadrature rule $\Ebb_\Ical$ is defined by summing over
all multi-indices in a prescribed multi-index set $\Ical \subset
\Nbb^{\ny}$, i.e.,
\begin{equation}
 \Ebb_\Ical \coloneqq \sum_{\ibm \in \Ical} \Delta^\ibm =
 \sum_{\ibm \in \Ical} \Delta_1^{i_1} \otimes \cdots \otimes
                       \Delta_{\ny}^{i_{\ny}}.
\end{equation}
An index set $\Ical$ must satisfy the standard admissibility requirement
\cite{gerstner2003dimension} for the above expression to be a convergent
quadrature approximation. This requirement is:  $\kbm - \ebm_j \in \Ical$
for all
$\kbm \in \Ical$,
$j\in\{\ell\mid 1\leq\ell\leq\ny,\, k_\ell >1\}$,
where $\ebm_j$ denotes the $j$th canonical unit
vector.
%satisfying $1 \leq j \leq \ny$.
%and $k_j > 1$ 

The collection of multi-indices that \emph{neighbor} a multi-index set
$\Ical $ is defined as
\begin{equation}
 \Ncal(\Ical) \coloneqq \{\ibm \in \Ical^c \mid \Ical \cup \{\ibm\}
                          \text{ is admissible}\}\subset \Nbb^{\ny},
\end{equation}
where $\Ical^c$ denotes the complement of the multi-index set $\Ical$ in
$\Nbb^{\ny}$, i.e.,
$\Ical^c \coloneqq \{\ibm \in \Nbb^{\ny} \mid \ibm \not\in \Ical\}$.
Following Refs.\
\cite{gerstner2003dimension, kouri2013trust, kouri2014inexact}, the truncation
error, which can be written as an infinite sum, can be approximated as
\begin{equation}
 \Ebb[h] - \Ebb_\Ical[h] = \sum_{\ibm \in \Ical^c}
 (\Delta_1^{i_1} \otimes \cdots \otimes \Delta_{\ny}^{i_{\ny}})
 [h]
 \approx \sum_{\ibm \in \Ncal(\Ical)}
 (\Delta_1^{i_1} \otimes \cdots \otimes \Delta_{\ny}^{i_{\ny}})
 [h].
\end{equation}

Finally, let $\dstochsp_\Ical \subset \dstochsp$ denote the quadrature nodes
associated with the multi-dimensional quadrature rule $\Ebb_\Ical$. If
nested one-dimensional quadrature rules are used, then
$\dstochsp_\Ical \subset \dstochsp_\Jcal$ for $\Ical,\,\Jcal$ multi-index sets
such that $\Ical \subset \Jcal$. This can lead to substantial savings since
evaluations of $h$ can be recycled as the sparse grid is refined.

\subsection{Projection-based model reduction}
\label{sec:aprxmod:rom}
In the construction of projection-based reduced-order models, the
distinction between the stochastic variables, $\dstochvar$, and
parameters, $\dparamvar$, is unimportant so, for notational brevity,
we combine them into a single vector,
\begin{equation}
 \dcombvar \coloneqq \begin{bmatrix} \dstochvar \\ \dparamvar \end{bmatrix}
\in \Rbb^{\nz},
\end{equation}
where $\nz\coloneqq\nzVerbose$.
\subsubsection{Minimum-residual primal reduced-order
model}\label{sec:minresprimal}
The projection-based reduced-order model construction begins with
the ansatz that the state vector lies in a low-dimensional subspace
\begin{equation} \label{eqn:rom-ansatz}
 \dpdesol \approx \rrob\dromsol,
\end{equation}
where $\rrob \in \Rbbstar^{\nU \times \kU}$ 
with $\kU \ll \nU$
denotes reduced basis matrix that spans the
low-dimensional subspace, $\dromsol \in \Rbb^{\kU}$ denotes the
generalized coordinates of the state $\dpdesol$ in the reduced
subspace, and $\Rbbstar^{m\times n}$ denotes the 
the set of full-column-rank $m\times n$ real-valued matrices
(i.e., the noncompact Stiefel manifold).  

The governing equations for the reduced-order model are obtained by
substituting the ansatz (\ref{eqn:rom-ansatz}) into the HDM governing
equations (\ref{eqn:sppde-disc}) and enforcing orthogonality of the
residual to  a test basis matrix
$\lrob:\Rbb^\nU\times\Rbb^{\nz}\rightarrow\Rbbstar^{\nU \times \kU}$, yielding
the following problem statement: Given $\dcombvar\in\Rbb^{\nz}$, compute 
the reduced primal solution $\dromsolSol$ satisfying
\begin{equation} \label{eqn:rom}
 \lrob(\rrob\dromsolSol,\,\dcombvar)^T\dpderes(\rrob\dromsolSol,\,\dcombvar) =
 \zerobold.
\end{equation}
We assume that for each $\dcombvar\in\Rbb^{\nz}$ there exists a unique reduced solution
$\dromsolSol = \dromsolSol(\dcombvar)$ satisfying \eqnref{eqn:rom} and
the implicit map $\dcombvar\mapsto \dromsolSol$ is continuously differentiable with respect to its arguments.

The reduced-order model in (\ref{eqn:rom}) is said to possess the
minimum-residual property \cite{legresley2006application,
                                carlberg2011gnat, zahr2016phd}
if it is equivalent to minimizing the HDM residual $\dpderes$ in
some metric over the low-dimensional trial space $\Range(\rrob)$,
with $\Range(\boldsymbol A)$ denoting the range of matrix $\boldsymbol A$.
\begin{definition}[Minimum-Residual Property] \label{def:minres}
A reduced-order model possesses the minimum-residual property if the solution
	$\dromsolSol$ satisfying Eq.~\eqref{eqn:rom} also
satisfies 
\begin{equation} \label{eqn:rom-minres-optform}
 \dromsolSol=
	\underset{\dromsol  \in \Rbb^\kU}{\arg\min}~
        {\frac{1}{2}\norm{\dpderes(\rrob\dromsol,\,\dcombvar)}_\Thetabold^2}
\end{equation}
given $\dcombvar\in\Rbb^{\nz}$ for some symmetric positive-definite matrix
$\Thetabold \in \Rbb^{\nU \times \nU}$.
\end{definition}
The first-order optimality conditions of the unconstrained optimization
problem \eqref{eqn:rom-minres-optform} provide a necessary condition for a
reduced-order model to possess the minimum-residual property, namely
\begin{equation} \label{eqn:lrob-minres}
 \lrob(\dpdesol,\,\dcombvar) =
 \Thetabold\pder{\dpderes}{\dpdesol}(\dpdesol,\,\dcombvar)\rrob,
\end{equation}
as this choice yields equivalence between the ROM governing equations
\eqref{eqn:rom} and the first-order optimality conditions of problem
\eqref{eqn:rom-minres-optform}.
Minimum-residual reduced-order models have the benefit of minimizing
the error bound on the QoI in (\ref{eqn:rom-qoi-err-indic}) and possess three
key properties defined in Proposition~\ref{prop:minres-monotone-interp}
\cite{zahr2016phd}:
\begin{inparaenum}[(1)]
 \item optimality,
 \item monotonicity, and
 \item interpolation.
\end{inparaenum}
\begin{prop} \label{prop:minres-monotone-interp}
 Let $(\rrob,\,\Thetabold)$ define a minimum-residual primal reduced-order model
 whose solution $\dromsolSol \in \Rbb^\kU$ for a given $\dcombvar \in \Rbb^{\nz}$,
 satisfies (\ref{eqn:rom-minres-optform}). Then,
 the following properties hold:
\begin{itemize}
 \item (Optimality) For any $\dpdesol\in\Range(\rrob)$,
\begin{equation}
 \norm{\dpderes(\rrob\dromsolSol,\,\dcombvar)}_\Thetabold \leq
 \norm{\dpderes(\dpdesol,\,\dcombvar)}_\Thetabold.
\end{equation}
 \item (Monotonicity) Let $\dromsolSol' \in \Rbb^\kU$ be the solution for a
	 given $\dcombvar$ corresponding to
        any projection-based reduced-order model with reduced
        basis $\rrob' \in \Rbb^{\nU\times\kU}$ such that
				$\Range(\rrob')\subseteq\Range(\rrob)$, then
\begin{equation}
 \norm{\dpderes(\rrob \dromsolSol ,\,\dcombvar)}_\Thetabold \leq
 \norm{\dpderes(\rrob'\dromsolSol',\,\dcombvar)}_\Thetabold.
\end{equation}
 \item (Interpolation) Let $\dpdesolSol$ satisfy
 $\dpderes(\dpdesolSol,\,\dcombvar) = \zerobold$. If $\dpdesolSol \in \Range(\rrob)$, then
\begin{equation} \label{eqn:minres-rom-interp-prop}
	\dpderes(\rrob\dromsolSol,\,\dcombvar) = \zerobold, \quad\text{with}\quad
 \dromsolSol = \rrob^+\dpdesolSol,
\end{equation}
where a superscript $+$ denotes the Moore--Penrose pseudo-inverse.
\end{itemize}
\begin{proof}
 Optimality follows directly from the assumption that the $\dromsolSol$
 is the global minimum of (\ref{eqn:rom-minres-optform}): since
 $\dromsolSol$ is the global minimum of the residual over the trial space,
 any other point in the trial space must have a larger residual. Monotonicity
 and interpolation follow directly from optimality.
\end{proof}
\end{prop}

The quantity of interest associated with any reduced-order model characterized by a unique mapping
$\dcombvar\mapsto\dromsolSol$ can be computed from the mapping
\begin{equation} \label{eqn:rom-qoi}
 \rqoirstr:\dcombvar \mapsto \qoi(\rrob\dromsolSol(\dcombvar),\,\dcombvar),
\end{equation}
which can be equipped with a residual-based error indicator
if the quantity of interest and residual are sufficiently regular,
as defined in Appendix~\ref{app:resbnd}, i.e.,
\begin{equation} \label{eqn:rom-qoi-err-indic}
 |\qoirstr(\dcombvar)-\rqoirstr(\dcombvar)| \leq
 \kappa \norm{\dpderes(\rrob\dromsolSol(\dcombvar),\,\dcombvar)}_\Thetabold,
\end{equation}
where $\kappa > 0$ is independent of $\dcombvar$,
$\Thetabold \in \Rbb^{\nU \times \nU}$ is a symmetric positive-definite
(SPD) matrix, and $\norm{\xbm}_\Thetabold^2 \coloneqq \xbm^T\Thetabold\xbm$ for any
$\xbm \in \Rbb^\nU$; see Appendix~\ref{app:resbnd}.

The problem statement associate with the adjoint of the reduced-order
model in (\ref{eqn:rom}) for a general test basis
$\lrob = \lrob(\dpdesol,\,\dcombvar)$ is the following: 
Given $\dcombvar\in\Rbb^{\nz}$ and
the reduced primal solution $\dromsolSol$ satisfying Eq.~\eqref{eqn:rom},
compute
the reduced adjoint solution $\dromadjSol \in \Rbb^\kU$ satisfying
\begin{equation} \label{eqn:rom-adj-gen}
 \left[\sum_{j=1}^{\nU}
       \dpderes_j\pder{\left(\lrob^T\ebm_j\right)}{\dpdesol}\rrob+
       \lrob^T\pder{\dpderes}{\dpdesol}\rrob\right]^T\dromadjSol =
 \rrob^T\pder{f}{\dpdesol}^T,
\end{equation}
where all quantities in Eq.~(\ref{eqn:rom-adj-gen}) are evaluated at
$(\rrob\dromsolSol,\,\dcombvar)$.
From the reduced primal--adjoint pair $(\dromsolSol,\,\dromadjSol)$ satisfying
Eqs.~\eqref{eqn:rom} and \eqref{eqn:rom-adj-gen},
the gradient of the reduced quantity of interest
can be computed as
\begin{equation} \label{eqn:rom-adj-qoi-gen}
 \nabla_\dparamvar\rqoirstr(\dcombvar) =
 \pder{\qoi}{\dparamvar} - 
 \dromadjSol^T
 \left[\sum_{j=1}^{\nU}\dpderes_j\pder{\left(\lrob^T\ebm_j\right)}{\dparamvar} +
       \lrob^T\pder{\dpderes}{\dparamvar}\right],
\end{equation}
where again all quantities are evaluated at
$(\rrob\dromsolSol,\,\dcombvar)$.

If the reduced primal solution $\dromsolSol$ is exact such that
$\dpderes(\rrob\dromsolSol(\dcombvar),\,\dcombvar) = \zerobold$
or the test basis $\lrob$ is independent of the primal state vector
$\dpdesol$, the reduced
adjoint equations (\ref{eqn:rom-adj-gen}) simplify to
\begin{equation} \label{eqn:rom-adj-exact}
 \left[\lrob^T\pder{\dpderes}{\dpdesol}\rrob\right]^T\dromadjSol =
 \rrob^T\pder{f}{\dpdesol}^T.
\end{equation}
and the gradient of the reduced quantity of interest simplifies to
\begin{equation}
	\nabla_\dparamvar\rqoirstr(\dcombvar) =
 \pder{\qoi}{\dparamvar}(\rrob\dromsolSol,\,\dcombvar) - 
 \dromadjSol^T\lrob(\rrob\dromsolSol,\,\dcombvar)^T\pder{\dpderes}{\dparamvar}(\rrob\dromsolSol,\,\dcombvar)=
	\dadjgrad(\lrob(\rrob\dromsolSol,\,\dcombvar)\dromadjSol,\,\rrob\dromsolSol,\,\dcombvar).
\end{equation}

The reduced adjoint equations (\ref{eqn:rom-adj-gen}) and the corresponding
expression to construct the gradient of the reduced quantity of interest
(\ref{eqn:rom-adj-qoi-gen}) can be challenging to implement, as computing the
derivative of the test basis entails computing the second derivative of the
residual. Furthermore, the reduced adjoint variable $\dromadjSol$ satisfying
Eq.~\eqref{eqn:rom-adj-gen} cannot necessarily be reconstructed in the full
space such that it provides an accurate approximation of the full-system
adjoint variable $\dadjvarSol$ satisfying Eq.~\eqref{eqn:adjres}, nor
does the reduced gradient  $\nabla_\dparamvar\rqoirstr$
accurately approximate its full-system counterpart 
$\nabla_\dparamvar\qoirstr$.
Therefore, instead of adhering to this
\emph{reduce-then-optimize} \cite{negri2015reduced} approach, we instead
construct a minimum-residual reduced-order model for the full-system
adjoint equations (\ref{eqn:adjres}). While this approach forfeits consistency
between the reduced quantity of interest $\rqoirstr$ and its derivative,
it gains the notion of residual minimization and its properties
(i.e., optimality, monotonicity, interpolation).

\subsubsection{Minimum-residual adjoint reduced-order model}
\label{sec:minresadjoint}
The minimum-residual adjoint construction begins with the ansatz that the
adjoint solution can be well approximated in a low-dimensional subspace.
For simplicity, we employ the same trial subspace to approximation the primal
and adjoint solutions
\begin{equation} \label{eqn:rom-adj-ansatz}
 \dadjvar \approx \rrob \dromadjmr,
\end{equation}
but this is not required. Reference \cite{zahr2016phd} takes the
adjoint trial subspace to be equal to the primal test subspace and shows this
to be
sufficient for the  minimum-residual adjoint and true reduced-order model
adjoint to agree whenever the reduced primal solution is exact.
The $\hat{\cdot}$ notation is used to distinguish the reduced-order adjoint
model from the true adjoint of the reduced system. 

We propose to compute the reduced adjoint variable by minimizing the
$\Thetabold^\dadjvar$-norm of the adjoint residual over the trial subspace
such that the associated problem statement becomes: 
Given $\dcombvar\in\Rbb^{\nz}$ and
the reduced primal solution $\dromsolSol$ satisfying Eq.~\eqref{eqn:rom},
compute the minimum-residual adjoint solution
$\dromadjmrSol\in\Rbb^{\kU}$ satisfying
\begin{equation} \label{eqn:rom-adj-optform1}
 \dromadjmrSol = 
 \underset{\etabold \in \Rbb^\kU}{\arg\min}
           ~\frac{1}{2}\norm{\dadjres(\rrob\etabold,\,\rrob\dromsolSol,\,\dcombvar)}_{\Thetabold^\dadjvar}^2.
\end{equation}
The first-order optimality conditions of (\ref{eqn:rom-adj-optform1}) can be
written as the linear system
\begin{equation} \label{eqn:romadj-minres-1st-opt}
 \left(\pder{\dpderes}{\dpdesol}^T\rrob\right)^T\Thetabold^\dadjvar
 \left(\pder{\dpderes}{\dpdesol}^T\rrob\right)\dromadjmrSol =
 \left(\pder{\dpderes}{\dpdesol}^T\rrob\right)^T\Thetabold^\dadjvar
 \pder{\qoi}{\dpdesol}^T,
\end{equation}
where all quantities in Eq.~(\ref{eqn:romadj-minres-1st-opt}) are evaluated at
$(\rrob\dromsolSol,\,\dcombvar)$.
 The minimum-residual adjoint reduced-order
model also possess similar concepts of optimality, monotonicity, and
interpolation as the primal \cite{zahr2016phd}.
\begin{prop} \label{prop:monotone-interp-adj-rom}
Let $(\rrob,\,\Thetabold^\dadjvar)$ define a
minimum-residual adjoint reduced-order model
whose solution	
$\dromadjmr\in \Rbb^\kU$ for a given $\dcombvar \in \Rbb^{\nz}$,
 satisfies (\ref{eqn:rom-adj-optform1}).
	Then the following properties hold:
\begin{itemize}
	\item (Optimality) For any $\dadjvar\in\Range(\rrob)$,
 \begin{equation}
  \norm{\dadjres(\rrob\dromadjmrSol,\,\rrob\dromsolSol,\,\dcombvar)}_{\Thetabold^\dadjvar}
  \leq
  \norm{\dadjres(\dadjvar,\,\rrob\dromsolSol,\,\dcombvar)}_{\Thetabold^\dadjvar}.
 \end{equation}
 \item (Monotonicity) Let $(\rrob',\,{\Thetabold^\dadjvar})$ define
       a minimum-residual adjoint reduced-order model such that\newline
       $\Range(\rrob') \subseteq \Range(\rrob)$, then
\begin{equation}
 \norm{\dadjres(\rrob\dromadjmrSol,\,\rrob\dromsolSol,\,\dcombvar)}_{\Thetabold^\dadjvar}
 \leq
 \norm{\dadjres(\rrob'\dromadjmrSol',\,\rrob\dromsolSol,\,\dcombvar)}_{\Thetabold^\dadjvar},
\end{equation}
where $\dromadjmrSol'$ satisfies
\begin{equation} \label{eqn:rom-adj-optform1Two}
 \dromadjmrSol' = 
 \underset{\etabold \in \Rbb^\kU}{\arg\min}
           ~\frac{1}{2}\norm{\dadjres(\rrob'\etabold,\,\rrob\dromsolSol,\,\dcombvar)}_{\Thetabold^\dadjvar}^2.
\end{equation}
 \item (Interpolation) If the adjoint solution $\dadjvarSol$ satisfying
	 Eq.~\eqref{eqn:adjres} also satisfies
 $\dadjvarSol \in \Range(\rrob)$, then
\begin{equation} \label{eqn:minres-romadj-interp-prop}
  \dadjres(\rrob\dromadjmrSol,\,\rrob\dromsolSol,\,\dcombvar)=\zerobold, \quad\text{with}\quad
  \dromadjmrSol=\rrob^+\dadjvarSol.
\end{equation}
\end{itemize}
\begin{proof}
 Optimality follows directly from problem \eqref{eqn:rom-adj-optform1}:
	because
 $\dromadjmrSol$ minimizes the adjoint residual over the trial space and
 the optimization problem is convex (i.e., it is a linear-least-squares
	problems characterized by a full-column-rank
 matrix), any other element of the trial space will yield a larger
 objective-function value. Monotonicity and interpolation follow directly from optimality.
\end{proof}
\end{prop}

Minimum-residual sensitivity reduced-order models and the corresponding
approximation of the gradient of the QoI can be constructed similarly.
We refer to Refs.~\cite{zahr2016phd,zahr2015progressive} for details.

We propose to approximate the gradient of the QoI using
the adjoint-based gradient operator $\dadjgrad$ defined in
\eqref{eqn:adjgradDef} as
\begin{equation} \label{eqn:minres-adj-grad}
  \nabla_\dparamvar \qoirstr(\dcombvar) \approx
  \widehat{\nabla_\dparamvar \rqoirstr}(\dcombvar) \coloneqq
  \dadjgrad(\rrob\dromsolSol,\,\rrob\dromadjmrSol,\,\dcombvar),
\end{equation}
where the minimum-residual primal $\dromsolSol$ and adjoint $\dromadjmrSol$ solutions satisfy
\eqref{eqn:rom-minres-optform} and \eqref{eqn:rom-adj-optform1}, respectively.
We emphasize that $\widehat{\nabla_\dparamvar \rqoirstr}\neq \nabla_\dparamvar
\rqoirstr$ in general.

From Appendix~\ref{app:resbnd}, any approximation of this form is equipped
with the following residual-based error bound
\begin{equation} \label{eqn:minres-adj-grad-bnd0}
  \norm{\nabla_\dparamvar \qoirstr(\dcombvar) -
        \dadjgrad(\rrob\dromsol,\,\rrob\dromadjmr,\,\dcombvar)} \leq
  \kappa\norm{\dpderes(\rrob\dromsol,\,\dcombvar)}_\Thetabold +
  \tau\norm{\dadjres(\rrob\dromadjmr,\,\rrob\dromsol,\,\dcombvar)}_{\Thetabold^\dadjvar},
\end{equation}
for some constants $\kappa,\,\tau > 0$, any
$\dromsol,\,\dromadjmr \in \Rbb^\kU$, and
any symmetric-positive-definite matrices $\Thetabold,\,\Thetabold^\dadjvar \in
\Rbb^{\nU \times \nU}$.
We justify the choice of employing the minimum-residual primal $\dromsolSol$ and adjoint
$\dromadjmrSol$ solutions in the gradient estimate \eqref{eqn:minres-adj-grad}
by noting that this choice
minimizes the gradient-error bound (\ref{eqn:minres-adj-grad-bnd0});
critically, this property is not generally shared with the exact
gradient of the reduced quantity of interest, i.e.,
$\nabla_\dparamvar\rqoirstr$.

\section{Approximation model and managed inexactness via trust-regions}
\label{sec:trammo}
We have now introduced two approximation techniques equipped with associated
error indicators: dimension-adaptive sparse grids, which enable efficient
quadrature in moderate-dimensional stochastic spaces, and projection-based
reduced-order models, which enable inexpensive PDE solutions. In this section,
we combine these into a single approximation model for the risk measure
\eqref{eqn:qoi-rstr-rm} associated with a stochastic, parametrized nonlinear
system \eqref{eqn:sppde-disc}. The level of the sparse grid and
dimension of the reduced-order model serve as two refinement mechanisms
that trade cost for improved approximation accuracy.
With this refinement-equipped approximation
model, we develop a trust-region-based approximation model management
method \cite{alexandrov1998trust} to efficiently solve the original
stochastic optimization problem (\ref{eqn:sppde-opt}). Critically, we
show that the resulting approach enables global convergence to the solution of
the original problem.

Section \ref{sec:trammo:tr} reviews details of the trust-region method
proposed in Ref.~\cite{kouri2014inexact}, which provides the foundation for
the proposed method. This approach allows for approximation models with
inexact objective and gradient evaluations, as well as flexible error
indicators. The remainder of the section formulates our approximation model
and error indicator in the context of the global convergence theory of
Ref.~\cite{kouri2014inexact} and introduces a refinement algorithm for the
approximation model to guarantee global convergence of our method.

\subsection{Trust-region method with inexact objective and gradient
            evaluations}
\label{sec:trammo:tr}
Let $\func{\obj}{\Rbb^\nmu}{\Rbb}$ be a smooth functional satisfying 
assumptions stated below. At trust-region iteration $k$ characterized by
trust-region center
$\dparamvar_k$, a trust-region method constructs a smooth approximation model
$\func{\aprxmod_k}{\Rbb^\nmu}{\Rbb}$ such that $\aprxmod_k\approx\obj$ within
the trust region
$\{\dparamvar \in \Rbb^\nmu \mid \norm{\dparamvar-\dparamvar_k}\leq\trrad_k\}$,
where $\trrad_k > 0$ denotes the trust-region radius. A trust-region method
computes a trial iterate
$\hat\dparamvar_k$ by approximately solving the trust-region
subproblem
\begin{equation}
 \optcona{\dparamvar \in \Rbb^\nmu}{\aprxmod_k(\dparamvar)}
         {\norm{\dparamvar-\dparamvar_k} \leq \trrad_k.}
\end{equation}
The computed trial iterate must satisfy the fraction of Cauchy decrease
condition \cite{conn2000trust}.

\subsubsection{The gradient condition}\label{sec:gradCond}
Following Refs.~\cite{heinkenschloss2002analysis, kouri2014inexact},
the gradient of the model $\aprxmod_k$ must sufficiently approximate that of
the true objective
function $\obj$ at the trust-region center $\dparamvar_k$, i.e.,
\begin{equation}\label{eq:accuracyRequirement}
 \norm{\nabla\aprxmod_k(\dparamvar_k)-\nabla\obj(\dparamvar_k)} \leq
 \xi \min\{\norm{\nabla\aprxmod_k(\dparamvar_k)},\,\trrad_k\},
\end{equation}
where $\xi > 0$ is independent of $k$. Suppose the model $\aprxmod_k$ is
equipped with the error bound
\begin{equation} \label{eqn:tr-graderr0}
 \norm{\nabla\aprxmod_k(\dparamvar_k)-\nabla\obj(\dparamvar_k)} \leq
 \xi \gradbnd_k(\dparamvar_k),
\end{equation}
where $\func{\gradbnd_k}{\Rbb^\nmu}{\Rbb}$ denotes an error indicator for the
gradient. Then requirement \eqref{eq:accuracyRequirement} can be
restated solely in terms of the error bound as
\begin{equation} \label{eqn:inexgrad-cond0}
	\gradbnd_k(\dparamvar_k) \leq
 \kappa_\gradbnd \min\{\norm{\nabla\aprxmod_k(\dparamvar_k)},\,\trrad_k\},
\end{equation}
with $\kappa_\gradbnd > 0$ a chosen constant, which we refer to as the
\textit{gradient condition}.

\subsubsection{The objective condition}\label{sec:objCond}
A trust-region method accepts the trial iterate $\hat\dparamvar_k$ if it
produces sufficient decrease in the objective function relative to the
decrease predicted by the model \cite{conn2000trust}. Traditionally,
the ratio of actual-to-predicted reduction is used to assess whether or not
the trial iterate is successful; this is
costly to compute in the present context because it requires evaluating
the true objective function at the trial iterate $\obj(\hat\dparamvar_k)$,
which in turn requires computing the exact expectation over stochastic space.
However, Ref.~\cite{kouri2014inexact} introduced a method for scrutinizing a
trust-region step that does not require evaluating the true objective
function yet nonetheless ensures global convergence. This technique introduces another
approximation $\func{\aprxmodobj_k}{\Rbb^\nmu}{\Rbb}$ for the true objective
function $\obj$ that must be equipped with the error bound
\begin{equation}\label{eq:theta_error_bound}
	|\obj(\dparamvar_k)-\obj(\hat\dparamvar_k)-(\aprxmodobj_k(\dparamvar_k)-
	\aprxmodobj_k(\hat\dparamvar_k))| \leq
 \sigma \objbnd_k(\hat\dparamvar_k)
\end{equation}
where $\sigma > 0$ is independent of $k$ and
$\func{\objbnd_k}{\Rbb^\nmu}{\Rbb}$ denotes an error indicator for the
decrease in the objective function from $\dparamvar_k$. To ensure global
convergence \cite{kouri2014inexact}, the following \textit{objective
condition} must hold
for a fixed $\omega \in (0,\,1)$:
\begin{equation} \label{eqn:inexobj-cond0}
 \objbnd_k^\omega(\hat\dparamvar_k) \leq
 \eta\min\{\aprxmod_k(\dparamvar_k)-\aprxmod_k(\hat\dparamvar_k),\,r_k\}.
\end{equation}
Here, $\eta < \min\{\eta_1,\,1-\eta_2\}$ and
$\{r_k\}_{k=1}^\infty \subset [0,\,\infty) \text{ such that}
\lim_{k\rightarrow\infty} r_k = 0$. Also,
$\eta_1$ and $\eta_2$ with
 $0 < \eta_1 < \eta_2 < 1$ are the thresholds used to determine
if the trial step should be accepted, which occurs if
\begin{equation}
\varrho_k \coloneqq \frac{\aprxmodobj_k(\dparamvar_k) -
                       \aprxmodobj_k(\hat\dparamvar_k)}
                      {\aprxmod_k(\dparamvar_k) - \aprxmod_k(\hat\dparamvar_k)}
                  \geq \eta_1,
\end{equation}
and if the trust-region radius should be increased, which occurs if
 $\varrho_k \geq \eta_2$. The complete algorithm, summarized in
Algorithm~\ref{alg:tr-iog-gen}, is globally convergent
\cite{kouri2014inexact} provided the following assumptions hold:
\begin{enumerate}
 \item $\obj$ is twice continuously differentiable and bounded below.
 \item $\aprxmodobj_k$ bounded below for all $k$.
 \item $\func{\aprxmod_k}{\Rbb^\nmu}{\Rbb}$ is twice continuously
       differentiable for all $k$.
 \item There exist $\kappa_1,\,\kappa_2 > 0$ such that, for all
       $\dparamvar \in \Rbb^\nmu$ and for all $k$,
 \begin{equation*}
  \norm{\nabla^2\obj(\dparamvar)} \leq \kappa_1 \qquad \text{and} \qquad
  \norm{\nabla^2\aprxmod_k(\dparamvar)} \leq \kappa_2.
 \end{equation*}
\end{enumerate}

\begin{algorithm}
 \caption{Trust region method with inexact objective evaluations
          \cite{kouri2014inexact}}
 \label{alg:tr-iog-gen}
 \begin{algorithmic}[1]
  \STATE \textbf{Initialization}: Given
 \begin{center}
  $\dparamvar_0$, $\trrad_0$,
  $0 < \gamma < 1$,
  %$\trrad_\text{max} > 0$,
  $0 < \eta_1 < \eta_2 < 1$, $0 < \kappa_\gradbnd$,
  \\
  $\omega \in (0,\,1)$,
  $\{r_k\}_{k=1}^\infty \subset [0, \infty)$ such that $r_k \rightarrow 0$
 \end{center}
	 \STATE\label{step:modelUpdate} \textbf{Model update}: Choose a model,
 $\aprxmod_k(\dparamvar)$, and gradient error bound,
 $\gradbnd_k(\dparamvar)$, such that
 \begin{equation*}
  \begin{aligned}
   \norm{\nabla \obj(\dparamvar_k) - \nabla \aprxmod_k(\dparamvar_k)} &\leq
   \xi \gradbnd_k(\dparamvar_k) \\
   \gradbnd_k(\dparamvar_k) &\leq \kappa_\gradbnd
                       \min\{\norm{\nabla \aprxmod_k(\dparamvar_k)},\,\trrad_k\}
  \end{aligned}
 \end{equation*}
 where $\xi > 0$ is an arbitrary constants
\vspace{1mm}
\STATE \textbf{Step computation}: Approximately solve the trust-region
 subproblem
 \begin{equation*}
  \underset{\dparamvar \in \Rbb^\nmu}{\min}~\aprxmod_k(\dparamvar) \qquad
  \text{subject to} \qquad \norm{\dparamvar-\dparamvar_k} \leq \trrad_k
 \end{equation*}
 for a candidate step $\hat\dparamvar_k$ that satisfies
 \begin{equation} %\label{eqn:fcd}
   \aprxmod_k(\dparamvar_k) - \aprxmod_k(\hat\dparamvar_k) \geq
   \kappa_s\norm{\nabla \aprxmod_k(\dparamvar_k)}
   \min\left\{\trrad_k,\,
              \frac{\norm{\nabla \aprxmod_k(\dparamvar_k)}}{\beta_k}\right\}
 \end{equation}
where $\kappa_s \in (0,\,1)$,
$\displaystyle{\beta_k \coloneqq 1 + \sup_{\dparamvar \in \Bcal_k}
               \norm{\nabla^2 \aprxmod_k(\dparamvar)}}$, and
$\Bcal_k \coloneqq \{\dparamvar \in \Rbb^{\nmu} \mid
                   \norm{\dparamvar-\dparamvar_k} \leq \trrad_k \}$
\vspace{1mm}
\STATE\label{step:actualToPredict} \textbf{Actual-to-predicted reduction}: Compute actual-to-predicted
 reduction ratio approximation
\begin{equation*}
 \varrho_k \coloneqq \frac{\aprxmodobj_k(\dparamvar_k)-\aprxmodobj_k(\hat\dparamvar_k)}
               {\aprxmod_k(\dparamvar_k)-\aprxmod_k(\hat\dparamvar_k)}
\end{equation*}
 where $\aprxmodobj_k(\dparamvar)$ and $\objbnd_k(\dparamvar)$ satisfy
\begin{equation}
 \begin{aligned}
  |\obj(\dparamvar_k)-\obj(\dparamvar)+
   \aprxmodobj_k(\dparamvar)-\aprxmodobj_k(\dparamvar_k)| &\leq
  \sigma \objbnd_k(\dparamvar), \qquad \dparamvar \in \Bcal_k \\
  \objbnd_k^\omega(\hat\dparamvar_k) &\leq
         \eta \min\{\aprxmod_k(\dparamvar_k)-\aprxmod_k(\hat\dparamvar_k),~r_k\}
 \end{aligned}
\end{equation}
where $\eta < \min\{\eta_1,\,1-\eta_2\}$ and $\sigma > 0$ is an arbitrary
constant
\vspace{1mm}
\STATE \textbf{Step acceptance}:
\begin{equation*}
 \textbf{if} \qquad \varrho_k \geq \eta_1 \qquad \textbf{then}
  \qquad \dparamvar_{k+1} = \hat\dparamvar_k \qquad \textbf{else}
  \qquad \dparamvar_{k+1} = \dparamvar_k \qquad \textbf{end if}
\end{equation*}
\STATE \textbf{Trust region update}:
\begin{equation*}
\begin{aligned}
\textbf{if} \qquad &\varrho_k \leq \eta_1 \qquad &\textbf{then} \qquad
&\trrad_{k+1} \in (0, \gamma \norm{\dparamvar-\hat\dparamvar_k}] \qquad
&\textbf{end if}\\
\textbf{if} \qquad &\varrho_k \in (\eta_1, \eta_2)\qquad &\textbf{then}\qquad
 &\trrad_{k+1} \in [\gamma \norm{\dparamvar-\hat\dparamvar_k}, \trrad_k] \qquad
 &\textbf{end if}\\
\textbf{if} \qquad &\varrho_k \geq \eta_2 \qquad &\textbf{then} \qquad
&\trrad_{k+1} \in [\trrad_k, +\infty) \qquad &\textbf{end if}\\
\end{aligned}
\end{equation*}
\end{algorithmic}
\end{algorithm}

The remainder of this section introduces approximation models
$\aprxmod_k$ and $\aprxmodobj_k$ and their associated error bounds
$\gradbnd_k$ and $\objbnd_k$, respectively, for the risk measure
associated with the quantity of interest \eqref{eqn:qoi-rstr-rm}
based on sparse grids and
reduced-order models. The section also develops an algorithm that ensures
global convergence of the method by adaptively refining the sparse grid and
reduced-order model such that the approximation models are sufficiently
accurate at trust-region centers according to the gradient condition
(\ref{eqn:inexgrad-cond0}) and the objective condition
(\ref{eqn:inexobj-cond0}).

\subsection{Approximation model based on sparse grids and reduced-order models}
\label{sec:trammo:aprxmod}
In this section, we combine the approximation techniques from
Section~\ref{sec:aprxmod} with the trust-region method described in the Section \ref{sec:trammo:tr}
to develop a globally convergent optimization method that leverages
and manages inexactness to solve the \textit{risk-neutral} optimization problem
\begin{equation} \label{eqn:sppde-opt-expect}
 \optunc{\dparamvar\in\Rbb^\nmu}
        {\obj(\dparamvar) \coloneqq
         \Ebb\left[\qoi(\dpdesolSol(\,\cdot\,,\,\dparamvar),\,
                        \cdot\,,\,\dparamvar)\right]},
\end{equation}
As indicated by
problem (\ref{eqn:sppde-opt-expect}), we return to the notation introduced in
Section~\ref{sec:prob} that distinguishes between stochastic variables
$\dstochvar$ and optimization parameters $\dparamvar$. % The development in this
%section is specific to the risk-neutral risk measure. however, the
%construction can be extended to general, smooth risk measures.

At the $k$th trust-region iteration, define the approximation models 
$\aprxmod_k$ and $\aprxmodobj_k$
discussed in Section~\ref{sec:trammo:tr} as
\begin{equation}\label{eq:approximationModels}
 \begin{aligned}
  \aprxmod_k(\dparamvar) &\coloneqq
  \Ebb_{\Ical_k}\left[f(\rrob_k\dromsol_k(\,\cdot\,,\,\dparamvar),\,
												\cdot\,,\,\dparamvar)\right],\quad\text{and} \\
  \aprxmodobj_k(\dparamvar) &\coloneqq
  \Ebb_{\Ical_k'}\left[f(\rrob_k'\dromsol_k'(\,\cdot\,,\,\dparamvar),\,
                         \cdot\,,\,\dparamvar)\right],
 \end{aligned}
\end{equation}
respectively,
%\DPK{The model $m_k$ should be quadratic if we are using Truncated CG.}
%\MJZ{Didn't you use a non-quadratic $m_k$ in \cite{kouri2013trust}
%     (Section 4.2) with TCG?} 
where $(\Ical_k,\,\rrob_k)$ and $(\Ical_k',\,\rrob_k')$ are
sparse-grid/reduced-basis pairs and, for any $\dstochvar \in \dstochsp$ and
$\dparamvar \in \Rbb^\nmu$, $\dromsol_k(\dstochvar,\,\dparamvar)$ denotes the
unique solution satisfying
$\lrob_k^T\dpderes(\rrob_k\dromsol_k,\,\dstochvar,\,\dparamvar) = \zerobold$ and
$\dromsol_k'(\dstochvar,\,\dparamvar)$ denotes the unique solution 
satisfying
$\lrob_k'^T\dpderes(\rrob_k'\dromsol_k',\,\dstochvar,\,\dparamvar) =
\zerobold$.
We assume that the projection uniquely determines the test basis
matrix $\lrob_k$ (resp.\ $\lrob_k'$) from the trial basis
matrix $\rrob_k$ (resp.\ $\rrob_k'$), e.g., Galerkin projection
leads to $\lrob_k = \rrob_k$, minimum-residual leads to
(\ref{eqn:lrob-minres}). Approximation models \eqref{eq:approximationModels}
comprise the sparse grid approximation of the expectation with
high-dimensional-model evaluations replaced by reduced-order-model evaluations. The choice of the sparse-grid/reduced-basis pairs will be
driven by accuracy requirements on the approximation models provided by the gradient condition (\ref{eqn:inexgrad-cond0}) and
the objective condition (\ref{eqn:inexobj-cond0});
we defer a complete
discussion to Section~\ref{sec:trammo:adapt}. From
Eq.~\eqref{eqn:rom-adj-qoi-gen}, the gradient of the approximation model is
\begin{equation}
 \nabla\aprxmod_k(\dparamvar) = \Ebb_{\Ical_k}
               \left[\nabla_\dparamvar\qoirstr_r(\,\cdot\,,\,\dparamvar)\right].
\end{equation}
%\begin{equation}
% \nabla\aprxmod_k(\dparamvar) =
%  \Ebb_{\Ical_k}\left[\dadjgrad(\rrob_k\dromadj_k(\,\cdot\,,\,\dparamvar),\,
%                                \rrob_k\dromsol_k(\,\cdot\,,\,\dparamvar),\,
%                                \cdot\,,\,\dparamvar)\right].
%\end{equation}
For this expression to hold,
the adjoint $\dromadj_k$ appearing in
$\nabla_\dparamvar\qoirstr_r(\,\cdot\,,\,\dparamvar)$ (see
Eq.~\eqref{eqn:rom-adj-qoi-gen}) must satisfy
the true reduced adjoint
equations (\ref{eqn:rom-adj-gen}); however, for the case of minimum-residual
primal reduced-order models, the reduced adjoint solution is difficult or
impractical to compute for the reasons discussed in Section \ref{sec:minresprimal}.
%to compute due to the need to evaluation Hessian terms
%that are large, sparse, third-order tensors and rarely available in practice.
Therefore, as described in Section \ref{sec:minresadjoint}, we replace the reduced adjoint solution with the minimum-residual
adjoint solution satisfying Eq.~(\ref{eqn:rom-adj-optform1}) and approximate the model
gradient according to
\begin{equation}
 \nabla\aprxmod_k(\dparamvar) \approx
 \widehat{\nabla\aprxmod}_k(\dparamvar) \coloneqq
  \Ebb_{\Ical_k}\left[\dadjgrad(
                            \rrob_k\hat\dromadj_k(\,\cdot\,,\,\dparamvar),\,
                            \rrob_k\dromsol_k(\,\cdot\,,\,\dparamvar),\,
                            \cdot\,,\,\dparamvar)\right].
\end{equation}
where $\hat\dromadj_k(\dstochvar,\,\dparamvar)$ uniquely
satisfies
\begin{equation} \label{eqn:rom-adj-optform1Iter}
 \hat\dromadj_k = 
 \underset{\etabold \in \Rbb^\kU}{\arg\min}
           ~\frac{1}{2}\norm{\dadjres(\rrob_k\etabold,\,\rrob_k\dromsol_k,\,\dstochvar,\,\dparamvar)}_{\Thetabold^\dadjvar}^2.
\end{equation}

\begin{rem}
The quadrature rule $\Ebb_{\Ical_k}$ will typically contain both positive and negative weights, which may lead to a bad optimization step. However, the trust-region algorithm guarantees convergence if the stated conditions are satisfied and therefore provides a natural mechanism for controlling such issues.
\end{rem}

\subsubsection{Error indicator for the gradient condition}
To properly embed the proposed approximation model in the globally convergent,
inexact trust-region framework of Section~\ref{sec:trammo:tr}, the
approximation models must be equipped with error bounds of the form
(\ref{eqn:tr-graderr0}) and (\ref{eq:theta_error_bound}). The error
indicator for the model gradient takes the form
\begin{equation} \label{eqn:gradbnd0}
 \gradbnd_k(\dparamvar) \coloneqq
                \beta_1\Ecal_{1,k}(\dparamvar) +
                \beta_3\Ecal_{3,k}(\dparamvar) +
                \beta_4\Ecal_{4,k}(\dparamvar)
\end{equation}
where $\beta_1,\,\beta_3,\,\beta_4 > 0$ are constants chosen to balance
the contribution of the three terms and the individual error terms account
for the primal reduced-order model error ($\Ecal_{1,k}$), the adjoint
reduced-order model error ($\Ecal_{3,k}$), and the quadrature truncation
error ($\Ecal_{4,k}$), i.e.,
\begin{equation} \label{eqn:gradbnd1}
 \begin{aligned}
  \Ecal_{1,k}(\dparamvar) &\coloneqq \Ebb_{\Ical_k\cup\Ncal(\Ical_k)}
                         \left[\norm{\dpderes(\rrob_k\dromsol_k(\,\cdot\,,\,
                                                                  \dparamvar),\,
                                              \cdot\,,\,\dparamvar)}\right] \\
  \Ecal_{3,k}(\dparamvar) &\coloneqq \Ebb_{\Ical_k\cup\Ncal(\Ical_k)}
                          \left[\norm{\dadjres(
                                \rrob_k\dromadj_k(\,\cdot\,,\,\dparamvar),\,
                                \rrob_k\dromsol_k(\,\cdot\,,\,\dparamvar),\,
                                \cdot\,,\,\dparamvar)}\right] \\
  \Ecal_{4,k}(\dparamvar) &\coloneqq \Ebb_{\Ncal(\Ical_k)}
                          \left[\norm{\dadjgrad(
                                \rrob_k\dromadj_k(\,\cdot\,,\,\dparamvar),\,
                                \rrob_k\dromsol_k(\,\cdot\,,\,\dparamvar),\,
                                \cdot\,,\,\dparamvar)}\right].
 \end{aligned}
\end{equation}
Each term can be efficiently computed since they only require primal and
adjoint reduced-order model solves over the collocation nodes of the
current sparse grid, $\Ical_k$, and its forward neighbors, $\Ncal(\Ical_k)$.

To demonstrate that the error indicator $\gradbnd_k(\dparamvar)$ as
prescribed in (\ref{eqn:gradbnd0}) provides an asymptotic error bound
for the model gradient, consider the error between the objective and
model gradient:
\begin{equation}
 \norm{\nabla\obj(\dparamvar)-\widehat{\nabla\aprxmod}_k(\dparamvar)} =
 \norm{\Ebb\left[\dadjgrad(\dadjvar(\,\cdot\,,\,\dparamvar),\,
                           \dpdesol(\,\cdot\,,\,\dparamvar),\,
                           \cdot\,,\,\dparamvar)\right] -
       \Ebb_{\Ical_k}\left[\dadjgrad(
                                \rrob_k\dromadj_k(\,\cdot\,,\,\dparamvar),\,
                                \rrob_k\dromsol_k(\,\cdot\,,\,\dparamvar),\,
                                \cdot\,,\,\dparamvar)\right]}.
\end{equation}
Employing $\Ebb = \Ebb_{\Ical_k} + \Ebb_{\Ical_k^c}$ and the triangle
inequality yields
\begin{equation}
 \begin{aligned}
  \norm{\nabla\obj(\dparamvar)-\widehat{\nabla\aprxmod}_k(\dparamvar)} &\leq
  \Ebb\left[\norm{\dadjgrad(\dadjvar(\,\cdot\,,\,\dparamvar),\,\dpdesol(\,\cdot\,,\,\dparamvar),\,
                            \cdot\,,\,\dparamvar) -
                  \dadjgrad(\rrob_k\dromadj_k(\,\cdot\,,\,\dparamvar),\,
                            \rrob_k\dromsol_k(\,\cdot\,,\,\dparamvar),\,
                            \cdot\,,\,\dparamvar)\right]} \\ &+
   \Ebb_{\Ical_k^c}\left[\norm{\dadjgrad(\rrob_k\dromadj_k(\,\cdot\,,\,\dparamvar),\,
                                         \rrob_k\dromsol_k(\,\cdot\,,\,\dparamvar),\,
                                         \cdot\,,\,\dparamvar)}\right].
 \end{aligned}
\end{equation}
From Proposition~\ref{prop:dual-bnd}, under the regularity and boundedness
assumptions in Assumptions~\ref{assume:res}--\ref{assume:qoi}, there exist
constant $\kappa',\,\tau' > 0$ such that
\begin{equation}
 \begin{aligned}
  \norm{\nabla\obj(\dparamvar)-\widehat{\nabla\aprxmod}_k(\dparamvar)} \leq
   \kappa'\Ebb\left[\norm{\dpderes(\rrob_k\dromsol_k(\,\cdot\,,\,\dparamvar),\,
                                   \cdot\,,\,\dparamvar)}\right] &+
   \tau'\Ebb\left[\norm{\dadjres(\rrob_k\dromadj_k(\,\cdot\,,\,\dparamvar),\,
                                 \rrob_k\dromsol_k(\,\cdot\,,\,\dparamvar),\,
                                 \cdot\,,\,\dparamvar)}\right] \\ &+
   \Ebb_{\Ical_k^c}\left[\norm{\dadjgrad(\rrob_k\dromadj_k(\,\cdot\,,\,\dparamvar),\,
                                         \rrob_k\dromsol_k(\,\cdot\,,\,\dparamvar),\,
                                         \cdot\,,\,\dparamvar)}\right].
 \end{aligned}
\end{equation}
Finally, because the integrals above cannot be computed exactly for general
integrands, we approximate them as $\Ebb_{\Ical_k^c} \approx \Ebb_{\Ncal(\Ical_k)}$ and
$\Ebb = \Ebb_{\Ical_k\cup\Ical_k^c} \approx \Ebb_{\Ical_k\cup\Ncal(\Ical_k)}$,
which yields
\begin{equation}
 \begin{aligned}
  \norm{\nabla\obj(\dparamvar)-\widehat{\nabla\aprxmod}_k(\dparamvar)} &\leq
   \kappa'\Ebb_{\Ical_k\cup\Ncal(\Ical_k)}
              \left[\norm{\dpderes(\rrob_k\dromsol_k(\,\cdot\,,\,\dparamvar),\,
                                   \cdot\,,\,\dparamvar)}\right] \\ &+
   \tau'\Ebb_{\Ical_k\cup\Ncal(\Ical_k)}
            \left[\norm{\dadjres(\rrob_k\dromadj_k(\,\cdot\,,\,\dparamvar),\,
                                 \rrob_k\dromsol_k(\,\cdot\,,\,\dparamvar),\,
                                 \cdot\,,\,\dparamvar)}\right] \\ &+
   \Ebb_{\Ncal(\Ical_k)}
              \left[\norm{\dadjgrad(\rrob_k\dromadj_k(\,\cdot\,,\,\dparamvar),\,
                                    \rrob_k\dromsol_k(\,\cdot\,,\,\dparamvar),\,
                                    \cdot\,,\,\dparamvar)}\right] + \epsilon,
 \end{aligned}
\end{equation}
where the last term $\epsilon$ is the truncation error associated
with using the forward neighborhood to approximate the expectations. An estimate of
$\epsilon$ could be obtained using additional layers of forward neighborhoods, e.g.,
$\Ncal(\Ical_k\cup\Ncal(\Ical_k))$; however, since $\epsilon \rightarrow 0$  as the
sparse grid is refined (provided the integrand for the gradient $\nabla\obj(\dparamvar)$ is sufficiently regular with respect to the uncertain inputs), we choose to neglect it.
From this bound and the definition of $\gradbnd_k$ in
(\ref{eqn:gradbnd0}), the required gradient error bound
in (\ref{eqn:tr-graderr0}) follows.
%Because we ignored the truncation error
%$\epsilon$, the gradient error bound is not rigorously established; however,
%our numerical experiments show the method is still convergent

\subsubsection{Error indicator for the objective condition}
Similarly, the error indicator for the decrease in the objective function is
defined as
\begin{equation} \label{eqn:objbnd0}
 \objbnd_k(\dparamvar) \coloneqq \alpha_1\left(\Ecal_{1,k}'(\dparamvar)+
                                           \Ecal_{1,k}'(\dparamvar_k)\right) +
                             \alpha_2\left(\Ecal_{2,k}'(\dparamvar)+
                                           \Ecal_{2,k}'(\dparamvar_k)\right)
\end{equation}
where $\alpha_1,\,\alpha_2 > 0$ are constants chosen to balance the
contributions of the two terms and the individual
error terms account for the error in the primal solution ($\Ecal_{1,k}'$)
and quadrature truncation error ($\Ecal_{2,k}'$), i.e.,
\begin{equation} \label{eqn:objbnd1}
 \begin{aligned}
  \Ecal_{1,k}'(\dparamvar) &\coloneqq \Ebb_{\Ical_k'\cup\Ncal(\Ical_k')}
                       \left[\norm{\dpderes(\rrob_k'\dromsol_k'(\,\cdot\,,\,
                                                                  \dparamvar),\,
                                            \cdot\,,\,\dparamvar)}\right] \\
  \Ecal_{2,k}'(\dparamvar) &\coloneqq \Ebb_{\Ncal(\Ical_k')}
                          \left[\left|\qoi(\rrob_k'\dromsol_k'(\,\cdot\,,\,
                                                                  \dparamvar),\,
                                           \cdot\,,\,\dparamvar)
                          \right|\right].
 \end{aligned}
\end{equation}
Both terms can be efficiently computed since they only require primal
reduced-order model solves over the collocation nodes of the current
sparse grid, $\Ical_k'$, and its forward neighbors, $\Ncal(\Ical_k')$.

To demonstrate that the error indicator $\objbnd_k(\dparamvar)$ as prescribed
in \eqref{eqn:objbnd0} provides
an asymptotic error bound for the objective decrease, we use the triangle
inequality to bound the error in the objective decrease from $\dparamvar_k$
to $\dparamvar$ by the error in the objective at both points, i.e.,
\begin{equation}
 \left|\obj(\dparamvar_k)-\obj(\dparamvar)-
	(\aprxmodobj_k(\dparamvar_k)-\aprxmodobj_k(\dparamvar))
       \right| \leq
 \left|\obj(\dparamvar_k)-\aprxmodobj_k(\dparamvar_k)\right| +
 \left|\obj(\dparamvar)-\aprxmodobj_k(\dparamvar)\right|.
\end{equation}
The objective error for any $\dparamvar \in \Rbb^\nmu$ takes the form
\begin{equation}
 \left|\obj(\dparamvar)-\aprxmodobj_k(\dparamvar)\right| =
 \left|\Ebb\left[\qoi(\dpdesolSol(\,\cdot\,,\,\dparamvar),\,\cdot\,,\,\dparamvar)\right] -
       \Ebb_{\Ical_k'}
       \left[\qoi(\rrob_k'\dromsol_k'(\,\cdot\,,\,\dparamvar),\,\cdot\,,\,\dparamvar)
       \right]
 \right|.
\end{equation}
Again, the simple relation $\Ebb = \Ebb_{\Ical_k'}+\Ebb_{\Ical_k'^c}$ and triangle
inequality bound this by two terms: one that accounts for the error in the
reduced-order model, and one that accounts for the quadrature truncation
error, i.e.,
\begin{equation}
 \left|\obj(\dparamvar)-\aprxmodobj_k(\dparamvar)\right| \leq
 \Ebb\left[\left|\qoi(\dpdesolSol(\,\cdot\,,\,\dparamvar),\,\cdot\,,\,\dparamvar) -
                 \qoi(\rrob_k'\dromsol_k'(\,\cdot\,,\,\dparamvar),\,\cdot\,,\,\dparamvar)
           \right|\right] +
 \Ebb_{\Ical_k'^c}
       \left[\left|
         \qoi(\rrob_k'\dromsol_k'(\,\cdot\,,\,\dparamvar),\,\cdot\,,\,\dparamvar)
       \right|\right].
\end{equation}
From Proposition~\ref{prop:prim-bnd}, under the regularity and boundedness
assumptions in Assumptions~\ref{assume:res}--\ref{assume:qoi}, there exists a
constant $\kappa'>0$ such that
\begin{equation}
 \left|\obj(\dparamvar)-\aprxmodobj_k(\dparamvar)\right| \leq
 \kappa'
 \Ebb\left[\norm{\dpderes(\rrob_k'\dromsol_k'(\,\cdot\,,\,\dparamvar),\,
                          \cdot\,,\,\dparamvar)}\right] +
 \Ebb_{\Ical_k'^c}
       \left[\left|
         \qoi(\rrob_k'\dromsol_k'(\,\cdot\,,\,\dparamvar),\,\cdot\,,\,\dparamvar)
       \right|\right].
\end{equation}
Finally, we approximate the integrals in the above expression 
as $\Ebb_{\Ical_k'^c} \approx \Ebb_{\Ncal(\Ical_k')}$ and
$\Ebb=\Ebb_{\Ical_k'\cup\Ical_k'^c} \approx
\Ebb_{\Ical_k'\cup\Ncal(\Ical_k')}$ to yield
\begin{equation}
 \left|\obj(\dparamvar)-\aprxmodobj_k(\dparamvar)\right| \leq
 \kappa'
 \Ebb_{\Ical_k'\cup\Ncal(\Ical_k')}
     \left[\norm{\dpderes(\rrob_k'\dromsol_k'(\,\cdot\,,\,\dparamvar),\,
                          \cdot\,,\,\dparamvar)}\right] +
 \Ebb_{\Ncal(\Ical_k')}
       \left[\left|
         \qoi(\rrob_k'\dromsol_k'(\,\cdot\,,\,\dparamvar),\,\cdot\,,\,\dparamvar)
       \right|\right] + \epsilon'.
\end{equation}
where the last term $\epsilon'$ is the truncation error associated
with using the forward neighborhood to approximate the expectations. Similar to the
gradient condition, an estimate of $\epsilon'$ could be obtained using additional layers
of forward neighborhoods, i.e., $\Ncal(\Ical_k\cup\Ncal(\Ical_k))$; however, since
$\epsilon' \rightarrow 0$ as the sparse grid is refined (provided the integrand for the objective function $\obj(\dparamvar)$ is sufficiently regular with respect to the uncertain inputs), we choose to neglect it.
The required objective decrease error bound in (\ref{eq:theta_error_bound})
follows from the above bound and the definition of $\objbnd_k$ in
(\ref{eqn:objbnd0})--(\ref{eqn:objbnd1}).
%\KTC{I think we are missing
%something here. We need to provide a bound now for
%$\left|\obj(\dparamvar)-\aprxmodobj_k(\dparamvar)\right|$. I suggest we
%clearly connect the constants in \eqref{eq:theta_error_bound} with the
%constants in the definition of $\theta_k$}
%\DPK{Should the Appendix be in the text body?  It seems like a majority of
%our contribution is there.} 

\subsection{Adaptive construction of sparse grid and reduced basis}
\label{sec:trammo:adapt}
With the approximation models based on sparse grid quadrature and
reduced-order model evaluations defined in the previous section,
along with the associated error indicators, this section defines 
algorithms to construct the sparse-grid/reduced-basis pairs
$(\Ical_k,\,\rrob_k)$ and $(\Ical_k',\,\rrob_k')$
such that the  gradient condition (\ref{eqn:inexgrad-cond0}) and
the objective condition (\ref{eqn:inexobj-cond0}), which are required for
convergence, are satisfied. 

The first algorithm constructs
$(\Ical_k,\,\rrob_k)$, which define model $\aprxmod_k(\dparamvar)$ according
to definitions \eqref{eq:approximationModels},
in Step \ref{step:modelUpdate} of Algorithm \ref{alg:tr-iog-gen}
to ensure satisfaction of the gradient condition
(\ref{eqn:inexgrad-cond0}); the
second algorithm constructs
$(\Ical_k',\,\rrob_k')$,
which define model $\aprxmodobj_k(\dparamvar_k)$ according
to definitions \eqref{eq:approximationModels},
in Step \ref{step:actualToPredict} of Algorithm \ref{alg:tr-iog-gen} (i.e.,
after a trial step $\hat\dparamvar_k$ is available)
to ensure satisfaction of the objective condition
(\ref{eqn:inexobj-cond0}). Both algorithms combine the concepts of
dimension-adaptive sparse grid construction \cite{gerstner2003dimension} and
greedy construction of a reduced basis
\cite{patera2007reduced,rozza2008reduced}, wherein (1) 
the outer loop refines the sparse grid 
based on the local truncation error in the neighboring index set, and (2) the
inner loop greedily samples over the current sparse grid to ensure the
reduced-order-model contributions to the error are sufficiently small.

%Recall from (\ref{eqn:gradbnd0}), the gradient error indicator is a weighted
%sum of three terms: the primal and adjoint reduced-order model error and the
%quadrature truncation error
%\begin{equation}
% \gradbnd_k(\dparamvar) =
% \beta_1\Ecal_{1,k}(\dparamvar) +
% \beta_3\Ecal_{3,k}(\dparamvar) +
% \beta_4\Ecal_{4,k}(\dparamvar).
%\end{equation}
\subsubsection{Refinement algorithm for satisfying the gradient
condition}\label{eq:refinementGradCond}
As described in Section \ref{sec:gradCond}, global convergence of the trust-region method in
Section~\ref{sec:trammo:tr} is predicated on gradient condition \eqref{eqn:inexgrad-cond0}.
A sufficient condition for the gradient condition to hold is that each
term of $\gradbnd_k(\dparamvar_k)$ as defined in Eq.~\eqref{eqn:gradbnd0} satisfies an appropriate fraction of the condition, i.e.,
\begin{equation} \label{eqn:gradbnd2}
 \begin{aligned}
  \Ecal_{1,k}(\dparamvar_k) &\leq \frac{\kappa_\gradbnd}{3\beta_1}
                    \min\{\norm{\nabla\aprxmod_k(\dparamvar_k)},\,\trrad_k\} \\
  \Ecal_{3,k}(\dparamvar_k) &\leq \frac{\kappa_\gradbnd}{3\beta_3}
                    \min\{\norm{\nabla\aprxmod_k(\dparamvar_k)},\,\trrad_k\} \\
  \Ecal_{4,k}(\dparamvar_k) &\leq \frac{\kappa_\gradbnd}{3\beta_4}
                    \min\{\norm{\nabla\aprxmod_k(\dparamvar_k)},\,\trrad_k\} \\
 \end{aligned}
\end{equation}
The purpose of the positive weights $\beta_1,\,\beta_3,\,\beta_4$, introduced
in the previous section, is to balance the individual contributions of the
error terms such that the uniform split above is justified. This decomposition
reduces the task of satisfying
gradient condition (\ref{eqn:inexgrad-cond0}) to ensuring each inequality
in (\ref{eqn:gradbnd2}) is satisfied. While the interplay between the
error terms in (\ref{eqn:gradbnd2}) and refinement of the sparse grid and reduced-order model
is highly coupled and fairly complex, the following
observations suggest an effective training strategy:
\begin{inparaenum}[(1)]
 \item for a fixed reduced-order model, $\Ecal_{4,k}$ decreases
       (possibly non-monotonically) as the sparse grid is refined, and
 \item for a fixed sparse grid, $\Ecal_{1,k}$ and $\Ecal_{3,k}$ decrease
       (possibly non-monotonically) as the reduced-order model is
       hierarchically refined.
\end{inparaenum}
Therefore, for a fixed reduced-order model, the proposed method adapts the sparse grid 
using the anisotropic dimension-adaptive approach proposed in Ref.~\cite{gerstner2003dimension}
to reduce the truncation error $\Ecal_{4,k}$ and, given the updated sparse
grid, the proposed method employs a variant of the classical greedy method
\cite{patera2007reduced,rozza2008reduced} to adapt the reduced-order
model to reduce the error terms $\Ecal_{1,k}$ and $\Ecal_{3,k}$. The algorithm
performs these steps iteratively until
inequalities
(\ref{eqn:gradbnd2}) are satisfied. %Before discussing the combined algorithm in
%detail, the individual components, namely dimension-adaptive construction of a
%sparse grid and greedy construction of a reduced-order model, are introduced.
%Therefore, the construction of the reduced-order model, for a fixed sparse
%grid, will proceed according to a variant of the classical greedy method
%\cite{patera2007reduced,rozza2008reduced} to target the error terms
%$\Ecal_{1,k}$ and $\Ecal_{3,k}$. For a fixed
%reduced-order model, the sparse grid will be adapted using the anisotropic
%dimension-adaptive approach \cite{gerstner2003dimension} to target
%$\Ecal_{4,k}$. These steps will be performed iteratively until the conditions
%in (\ref{eqn:gradbnd2}) are met. Before discussing the combined algorithm in
%detail, the individual components, namely dimension-adaptive construction of a
%sparse grid and greedy construction of a reduced-order model, are introduced.

\textbf{Sparse-grid construction}
The proposed construction of $\Ical_k$ mimics the dimension-adaptive algorithm
introduced by Gerstner and Griebel
\cite{gerstner2003dimension} for constructing a goal-oriented, anisotropic
sparse grid. This approach approximates the truncation error associated with sparse
grid $\Ical_k$ using only the neighbors $\Ncal(\Ical_k)$.
If this truncation-error approximation is larger than a specified tolerance,
the multi-index in the set of neighbors that contributes most to the error
is added to the index set, i.e., $\Ical_k \leftarrow \Ical_k \cup \{\ibm^*\}$
with
\begin{equation} \label{eqn:max-idxset-truncerr}
 \ibm^* = \underset{\ibm\in\Ncal(\Ical_k)}{\arg\max}~\left|\Delta^\ibm[g]\right|
\end{equation}
and the integrand is $\func{g}{\Rbb^\ny}{\Rbb}$. In the context of the
proposed SG-ROM approximation, the dimension-adaptive algorithm is applied
to the integrand that appears in the gradient truncation error
$\Ecal_{4,k}(\dparamvar)$, i.e., for a fixed reduced basis matrix $\rrob_k$
and given $\dparamvar \in \Rbb^\nmu$, we set the integrand to
%\begin{equation*}
%\norm{\dadjgrad(\rrob_k\dromadjmr_k(\,\cdot\,,\,\dparamvar),\,
%                \rrob_k\dromsol_k(\,\cdot\,,\,\dparamvar),\,
%                \cdot\,,\,\dparamvar)},
%\end{equation*}
\begin{equation*}
	g:\dstochvar\mapsto\norm{\dadjgrad(\rrob_k\dromadjmr_k(\,\dstochvar\,,\,\dparamvar),\,
                \rrob_k\dromsol_k(\,\dstochvar\,,\,\dparamvar),\,
                \dstochvar\,,\,\dparamvar)}.
\end{equation*}
While convergence is not necessarily monotonic, this term approaches
zero as $\Ical_k \rightarrow \Nbb^\ny$.

\textbf{Reduced-basis construction}
The construction of the reduced basis follows the well-studied greedy
algorithm \cite{patera2007reduced, rozza2008reduced}. The original greedy
algorithm aims to improve the parametric robustness of a reduced-order model by
iteratively enriching the reduced basis with snapshots of the high-dimensional
model at the point in parameter space where the reduced-order-model error is
largest. Pragmatically, this is executed by evaluating the reduced-order model
and an inexpensive error indicator at a (possibly large) set of candidate
points in parameter space, and subsequently performing a direct search for the
maximum error-indicator value over this candidate set.  Ref.\
\cite{chen2014weighted} developed a weighted variant of the greedy algorithm
for stochastic problems with non-uniform probability distributions. This
approach uses the probability density $\probdens(\dstochvar)$ to weight the
error indicator while training the reduced-order model over the stochastic space
$\dstochsp$. This is sensible, as stochastic-space regions with larger
probability density yield larger contributions to the expected
reduced-order-model error. Ref.\ \cite{chen2014weighted} also coupled the
weighted greedy algorithm with sparse grids by using the quadrature nodes as
the candidate set; this is motivated by the observation that the reduced-order
model is queried only at these points in stochastic space.

This work adopts a similar weighted greedy algorithm to train
the reduced-order model over the quadrature nodes and neighbors
$\dstochsp_{\Ical\cup\Ncal(\Ical)}$ assuming fixed parameters $\dparamvar$ and
sparse grid $\Ical$. Because the gradient condition is required to hold only at
the trust-region center, the proposed method performs training solely in stochastic space,
with $\dstochsp_{\Ical\cup\Ncal(\Ical)}$ comprising the candidate set, for
$\dparamvar = \dparamvar_k$ fixed. Unlike traditional greedy methods, the
proposed method constructs a reduced basis that accurately represents
both the primal \emph{and} adjoint states over the candidate set. This
is required because the same reduced basis $\rrob$ is employed for the primal
\eqref{eqn:rom} and
adjoint \eqref{eqn:rom-adj-optform1} reduced-order models, and thus the greedy
algorithm is responsible for reducing both the
primal $\Ecal_{1,k}$ and adjoint $\Ecal_{3,k}$ error terms that
arise in the gradient error indicator $\gradbnd_k(\dparamvar)$
defined in
(\ref{eqn:gradbnd0}). The approach
achieves this
by adding both primal and adjoint snapshots to the reduced basis at each greedy iteration. 

From the form of the gradient
error indicator $\gradbnd_k$ in
(\ref{eqn:gradbnd0}), the primal error indicator is taken
as the primal residual norm weighted by the density, i.e., the integrand in
$\Ecal_{1,k}$ evaluated at $\dparamvar = \dparamvar_k$
\begin{equation} \label{eqn:prim-res-err}
 \probdens(\dstochvar)
 \norm{\dpderes(\rrob_k\dromsol_k(\dstochvar,\,\dparamvar_k),\,
                \dstochvar,\,\dparamvar_k)}_\Thetabold,
\end{equation}
and the adjoint error indicator is taken as the adjoint residual norm weighted by
the density,
i.e., the integrand in $\Ecal_{3,k}$ evaluated at $\dparamvar = \dparamvar_k$
\begin{equation} \label{eqn:dual-res-err}
  \probdens(\dstochvar)
  \norm{\dadjres(\rrob_k\dromadjmr_k(\dstochvar,\,\dparamvar_k),\,
                 \rrob_k\dromsol_k(\dstochvar,\,\dparamvar_k),\,
                 \dstochvar,\,\dparamvar_k)}_{\Thetabold^\dadjvar}.
\end{equation}
%for metrics $\Thetabold,\,\Thetabold^\dadjvar \succ 0$.

\textbf{Refinement algorithm}
With these error indicators, the dimension-adaptive algorithm to construct
the sparse grid $\Ical_k$ and reduced basis $\rrob_k$ at the $k$th
trust-region iteration takes form. The algorithm initializes these entities
to their values from the previous trust-region iteration as
$\Ical_k \leftarrow \Ical_{k-1}'$ and $\rrob_k \leftarrow \rrob_{k-1}'$.
Then, the algorithm expands the sparse grid using the quadrature nodes in
the neighboring index set that contribute most significantly to the
truncation error, i.e.,
\begin{equation} \label{eqn:dimadapt-grad}
 \begin{aligned}
  \ibm^* &= \underset{\ibm\in\Ncal(\Ical_k)}{\arg\max}~~\left|\Delta^\ibm[
                \norm{\dadjgrad(\rrob_k\dromadjmr_k(\,\cdot\,,\,\dparamvar_k),\,
                                \rrob_k\dromsol_k(\,\cdot\,,\,\dparamvar_k),\,
                                \cdot\,,\,\dparamvar_k)}]\right| \\
  \Ical_k &\leftarrow \Ical_k \cup \ibm^*.
 \end{aligned}
\end{equation}
With the updated sparse grid, the greedy algorithm uses the primal error
indicator
\begin{equation} \label{eqn:greedy-grad-prim}
 \begin{aligned}
  \dstochvar^* &=\underset{\dstochvar \in \dstochsp_{\Ical_k\cup\Ncal(\Ical_k)}}
                          {\arg\max}~
          \probdens(\dstochvar)
          \norm{\dpderes(\rrob_k\dromadjmr_k(\,\dstochvar\,,\,\dparamvar_k),\,
                         \rrob_k\dromsol_k(\,\dstochvar\,,\,\dparamvar_k),\,
                         \dstochvar\,,\,\dparamvar_k)}_\Thetabold \\
  \rrob_k &\leftarrow
  \begin{bmatrix}
   \rrob_k & \dpdesolSol(\dstochvar^*,\,\dparamvar_k) &
             \dadjvarSol(\dstochvar^*,\,\dparamvar_k)
  \end{bmatrix},
 \end{aligned}
\end{equation}
to reduce $\Ecal_{1,k}$ error below its required threshold in \eqref{eqn:gradbnd2}.
If a minimum-residual reduced-order model is employed, the algorithm is
guaranteed to terminate due to the monotonicity property in
Proposition~\ref{prop:minres-monotone-interp}. In the limiting case
where snapshots have been added for each
$\dstochvar \in \dstochsp_{\Ical\cup\Ncal(\Ical)}$, the primal reduced-order
model will be exact for each
$\dstochvar \in \dstochsp_{\Ical\cup\Ncal(\Ical)}$ and thus
$\Ecal_{1,k}=0$. Furthermore, the minimum-residual
property guarantees $\Ecal_{1,k}$ will decrease monotonically from the
monotonicity property in Proposition~\ref{prop:minres-monotone-interp}.

After 
$\Ecal_{1,k}$ has been reduced below its required threshold,
the algorithm employs another greedy method with the adjoint error indicator
\begin{equation} \label{eqn:greedy-grad-adj}
 \begin{aligned}
  \dstochvar^* &=\underset{\dstochvar \in \dstochsp_{\Ical_k\cup\Ncal(\Ical_k)}}
                          {\arg\max}~
          \probdens(\dstochvar)
          \norm{\dadjres(\rrob_k\dromadjmr_k(\,\dstochvar\,,\,\dparamvar_k),\,
                         \rrob_k\dromsol_k(\,\dstochvar\,,\,\dparamvar_k),\,
                         \dstochvar\,,\,\dparamvar_k)}_{\Thetabold^\lambdabold} \\
 \rrob_k &\leftarrow
 \begin{bmatrix}
  \rrob_k & \dpdesolSol(\dstochvar^*,\,\dparamvar_k) &
            \dadjvarSol(\dstochvar^*,\,\dparamvar_k)
 \end{bmatrix},
 \end{aligned}
\end{equation}
to reduce $\Ecal_{3,k}$ below its required threshold in
\eqref{eqn:gradbnd2}. Even if a minimum-residual adjoint
reduced-order model is used, $\Ecal_{3,k}$ is not guaranteed to decrease
monotonically because modification of the trial basis matrix $\rrob_k$ alters
the linearization point defining the adjoint residual. This can be seen
from Proposition~\ref{prop:monotone-interp-adj-rom} where the three properties
hold provided the primal state is fixed. However, in the limit where
all quadrature nodes in $\dstochsp_{\Ical\cup\Ncal(\Ical)}$ have been sampled,
the primal reduced-order model and adjoint reduced-order models will be
exact and the algorithm will terminate with $\Ecal_{3,k} = 0$.

Once the reduced basis is constructed such that $\Ecal_{1,k}$ and
$\Ecal_{3,k}$ satisfy their required bounds in (\ref{eqn:gradbnd2}),
the dimension-adaptive greedy algorithm terminates if the truncation
error indicator $\Ecal_{4,k}$ also satisfies its bound in (\ref{eqn:gradbnd2}). Otherwise, another iteration of the algorithm is performed, i.e., the
algorithm expands the sparse grid according to (\ref{eqn:dimadapt-grad})
and enriches the reduced basis according to
(\ref{eqn:greedy-grad-prim})--(\ref{eqn:greedy-grad-adj}).
Because the truncation error tends to zero as the sparse grid is refined,
this iteration is guaranteed to terminate provided the regularity conditions
on the state and adjoint described in 
Assumptions~\ref{assume:res}--\ref{assume:qoi} of
Appendix~\ref{app:resbnd} are met. Assuming the
algorithm terminates, the resulting sparse grid $\Ical_k$ and reduced basis
$\rrob_k$ are guaranteed to satisfy the gradient conditions in
(\ref{eqn:gradbnd2}), thereby ensuring the gradient condition
\eqref{eqn:inexgrad-cond0} required for global convergence is satisfied.

Algorithm~\ref{alg:refine-grad} summarizes the refinement algorithm
that applies the nested iteration described
above to satisfy the individual gradient
conditions (\ref{eqn:gradbnd2}). At trust-region iteration $k = 0$,
the algorithm initializes the sparse grid as the uniform level-one sparse grid
$\Ical = \{(1,\,\dots,\,1)\}$, which corresponds to a single quadrature node
$\Xibold_\Ical = \{\zerobold\}$\footnote{We assume the origin is the center of
$\dstochsp$ and the 1D quadrature rules have $\zerobold$ as the level 1
quadrature point, which is the case for several of the most common
nested 1D quadrature rules, including the Clenshaw-Curtis rules used
in this work.}. The algorithm constructs the reduced basis from
the primal and adjoint snapshot at this single quadrature node
computed at the first trust-region center, i.e., the snapshots are
$\dpdesolSol(\zerobold,\,\dparamvar_0)$
and $\dadjvarSol(\zerobold,\,\dparamvar_0)$. Then, the approach may refine the sparse grid and
reduced basis to yield $\Ical_0$
and $\rrob_0$, respectively, by applying
Algorithm~\ref{alg:refine-grad} with inputs $\Ical_{-1}' = \{(1,\,\dots,\,1)\}$
and $\rrob_{-1}' = [\dpdesolSol(\zerobold,\,\dparamvar_0),\,
                    \dadjvarSol(\zerobold,\,\dparamvar_0)]$.
For all subsequent trust-region iterations, the sparse grid and reduced
basis are initialized from their values at the previous trust-region iteration.
In addition to providing a natural way to initialize the dimension-adaptive
greedy algorithm, this choice has the following benefit: it expands the
sparse grid and enriches the reduced basis
only if the choices $\Ical_k = \Ical_{k-1}'$ and
$\rrob_k = \rrob_{k-1}'$ are not sufficient to guarantee convergence.
%\DPK{The discussion of the quadrature rule seems accurate only if
%the origin is the center of $\dstochsp$ and the 1D quadrature rules
%have 0 as the level 1 quadrature point.}

\begin{algorithm}
 %\caption{Refine reduced basis and sparse grid for gradient condition}
 \caption{Refinement algorithm for satisfying the gradient condition}
 \label{alg:refine-grad}
 \begin{algorithmic}[1]
 %\STATEx $$\Ical_k,\,\rrob_k =
 %          \texttt{nested-refine-grad}
 %          (\Ical_{k-1},\,\rrob_{k-1},\,\dparamvar_k,\,\trrad)$$
 \REQUIRE
  $\Ical_{k-1}'$, $\rrob_{k-1}'$, $\dparamvar_k$, $\trrad_k$, $\beta_1 > 0$,
  $\beta_3 > 0$, $\beta_4 > 0$, $\kappa_\gradbnd > 0$
 \ENSURE $\Ical_k$, $\rrob_k$
 \STATE \textbf{Set}:
        $\Ical_k \leftarrow \Ical_{k-1}'$, $\rrob_k \leftarrow \rrob_{k-1}'$
\WHILE{$\displaystyle{
        \Ecal_{4,k}(\dparamvar_k) >
        \frac{\kappa_\gradbnd}{3\beta_4}
        \min\left\{\norm{\Ebb_{\Ical_k}\left[\dadjgrad(
                                \rrob_k\dromadjmr_k(\,\cdot\,,\,\dparamvar_k),\,
                                \rrob_k\dromsol_k(\,\cdot\,,\,\dparamvar_k),\,
                                \cdot\,,\,\dparamvar_k)
                                                  \right]},\,\trrad_k\right\}}$}
 \vspace{1.5mm}
 \STATE \textbf{Refine index set}: Add index set with largest contribution to
         truncation error:
\begin{equation*}
 \Ical_k \leftarrow \Ical_k \cup \{\ibm^*\} \qquad\text{with}\qquad
 \ibm^* = \underset{\ibm \in \Ncal(\Ical_k)}{\arg\max}
          \left|\Delta^\ibm \left[\norm{[\dadjgrad(
                                \rrob_k\dromadjmr_k(\,\cdot\,,\,\dparamvar_k),\,
                                \rrob_k\dromsol_k(\,\cdot\,,\,\dparamvar_k),\,
                                \cdot\,,\,\dparamvar_k)}\right]\right|.
\end{equation*}
 \WHILE{$\displaystyle{
         \Ecal_{1,k}(\dparamvar_k) >
         \frac{\kappa_\gradbnd}{3\beta_1}
         \min\left\{\norm{\Ebb_{\Ical_k}\left[\dadjgrad(
                                \rrob_k\dromadjmr_k(\,\cdot\,,\,\dparamvar_k),\,
                                \rrob_k\dromsol_k(\,\cdot\,,\,\dparamvar_k),\,
                                \cdot\,,\,\dparamvar_k)
                                                  \right]},\,\trrad_k\right\}}$}
 \vspace{1.5mm}
 \STATE \textbf{Evaluate primal error indicator}: Greedily select
				$\ybm \in \Xibold_{\ibm^*}$ with largest error indicator:
 \begin{equation*}
  \dstochvar^* = \underset{\dstochvar \in \dstochsp_{\ibm^*}}{\arg\max}~
          \probdens(\dstochvar)
          \norm{\dpderes(\rrob_k\dromsol_k(\,\cdot\,,\,\dparamvar_k),\,
                         \cdot\,,\,\dparamvar_k)}_\Theta.
 \end{equation*}
 \STATE \textbf{Reduced-order model construction}: Update reduced basis with
         new snapshots:
 \begin{equation*}
  \rrob_k \leftarrow
  \begin{bmatrix}
   \rrob_k & \dpdesolSol(\dstochvar^*,\,\dparamvar_k) &
             \dadjvarSol(\dstochvar^*,\,\dparamvar_k)
  \end{bmatrix}
 \end{equation*}.
 \ENDWHILE
 \vspace{2mm}
 \WHILE{$\displaystyle{
         \Ecal_{3,k}(\dparamvar_k) >
         \frac{\kappa_\gradbnd}{3\beta_3}
         \min\left\{\norm{\Ebb_{\Ical_k}\left[\dadjgrad(
                                \rrob_k\dromadjmr_k(\,\cdot\,,\,\dparamvar_k),\,
                                \rrob_k\dromsol_k(\,\cdot\,,\,\dparamvar_k),\,
                                \cdot\,,\,\dparamvar_k)
                                                  \right]},\,\trrad_k\right\}}$}
 \vspace{1.5mm}
 \STATE \textbf{Evaluate dual error indicator}: Greedily select
        $\ybm \in \Xibold_{\ibm^*}$ with largest error:
 \begin{equation*}
  \dstochvar^* = \underset{\dstochvar \in \dstochsp_{\Ical_k\cup\Ncal(\Ical_k)}}
                          {\arg\max}~
          \probdens(\dstochvar)
          \norm{\dadjres(\rrob_k\dromadjmr_k(\,\cdot\,,\,\dparamvar_k),\,
                         \rrob_k\dromsol_k(\,\cdot\,,\,\dparamvar_k),\,
                         \cdot\,,\,\dparamvar_k)}_{\Theta^\lambdabold}
 \end{equation*}
 \STATE \textbf{Reduced-order model construction}: Update reduced basis with
         new snapshots:
 \begin{equation*}
  \rrob_k \leftarrow
  \begin{bmatrix}
   \rrob_k & \dpdesolSol(\dstochvar^*,\,\dparamvar_k) &
             \dadjvarSol(\dstochvar^*,\,\dparamvar_k)
  \end{bmatrix}
 \end{equation*}
 \ENDWHILE
\ENDWHILE
\end{algorithmic}
\end{algorithm}

This completes the discussion of the refinement algorithm used to construct $\Ical_k$
and $\rrob_k$ such that the gradient condition (\ref{eqn:inexgrad-cond0}) is
satisfied. We now turn attention to 
constructing $\Ical_k'$ and $\rrob_k'$
such that the objective condition (\ref{eqn:inexobj-cond0}) holds; this
enables assessment of the trust-region step without requiring queries to
$\obj(\dparamvar)$. 

\subsubsection{Refinement algorithm for satisfying the objective condition}

As described in Section \ref{sec:objCond},  
satisfying the objective condition (\ref{eqn:inexobj-cond0})
is also required to preserve global convergence of the trust-region method when
the approximate objective function $\aprxmodobj_k(\dparamvar)$ is used in
place of true objective function $\obj(\dparamvar)$ in the
computation of the actual-to-prediction reduction $\varrho_k$.  
Analogously to
inequalities \eqref{eqn:gradbnd2}, a sufficient condition for the objective condition to hold is that each
term of the error indicator $\objbnd_k(\dparamvar_k)$ as defined in
Eq.~\eqref{eqn:objbnd0} satisfies an appropriate fraction of the condition,
i.e.,
\begin{equation} \label{eqn:inexobj-cond-indiv}
 \begin{aligned}
  \Ecal_{1,k}'(\dparamvar_k) +
  \Ecal_{1,k}'(\hat\dparamvar_k)
  &\leq \frac{1}{2\alpha_1} \left(\eta
  \min\{\aprxmod_k(\dparamvar_k)-
        \aprxmod_k(\hat\dparamvar_k),\,r_k\}\right)^{1/\omega} \\
  \Ecal_{2,k}'(\dparamvar_k) +
  \Ecal_{2,k}'(\hat\dparamvar_k)
  &\leq \frac{1}{2\alpha_2} \left(\eta
  \min\{\aprxmod_k(\dparamvar_k)-
        \aprxmod_k(\hat\dparamvar_k),\,r_k\}\right)^{1/\omega}.
 \end{aligned}
\end{equation}
This decomposition reduces the monolithic task of satisfying the
objective condition \eqref{eqn:inexobj-cond0} to ensuring each
inequality in \eqref{eqn:inexobj-cond-indiv} is satisfied.
Analogously to the approach described in Section \ref{eq:refinementGradCond}
to construct $\Ical_k$ and $\rrob_k$, we propose to employ a weighted greedy
algorithm to enforce inequalities \eqref{eqn:inexobj-cond-indiv}.

While proposed dimension-adaptive greedy algorithm to construct
$(\Ical_k',\,\rrob_k')$ for the objective decrease condition will
be very similar to that used to construct $(\Ical_k,\,\rrob_k)$, there
will be two critical differences. First, the error terms in
(\ref{eqn:inexobj-cond-indiv}) involve \emph{two} points in parameters
space: the trust-region center $\dparamvar_k$ and the candidate step
$\hat\dparamvar_k$. In contrast, the error terms in the gradient-condition 
inequalities \eqref{eqn:gradbnd2} impose requirements only at the trust-region
center $\dparamvar_k$. This has implications for both the dimension-adaptive
sparse grid construction and greedy method. Second, comparing definitions
\eqref{eqn:gradbnd1} and \eqref{eqn:objbnd1} reveals that inequalities
(\ref{eqn:inexobj-cond-indiv}) impose requirements only on the primal
reduced-order model accuracy (through $\Ecal_{1,k}'$) and truncation
error (through $\Ecal_{2,k}'$), whereas the gradient-condition inequalities
\eqref{eqn:gradbnd2} also placed requirements on the
adjoint accuracy (through $\Ecal_{3,k}$). This implies only \emph{primal}
snapshots are required during the greedy construction of the reduced-order
model. However, we choose to also include adjoint snapshots because
the candidate $\hat\dparamvar_k$ will become the trust-region center at
the next iteration $\dparamvar_{k+1}$ if the
trust-region step is successful, and adjoint accuracy is required at
trust-region centers to satisfy the gradient condition.

\textbf{Sparse-grid construction}
 For a given $(\Ical_k',\,\rrob_k')$, if the truncation error
conditions in inequalities (\ref{eqn:inexobj-cond-indiv}), i.e., the requirements on
$\Ecal_{2,k}'$, are not satisfied, the algorithm updates the sparse grid according to
$\Ical_k' \leftarrow \Ical_k' \cup \{\ibm^*\}$, where
\begin{equation}
 \ibm^* = \underset{\ibm\in\Ncal(\Ical)}{\arg\max}
 \left(\max\left\{
 \left|\Delta^\ibm\left[
  \qoi(\rrob_k'\dromsol_k'(\,\cdot\,,\,\dparamvar_k),\,\cdot\,,\,\dparamvar_k)
  \right]\right|,\,
 \left|\Delta^\ibm\left[
  \qoi(\rrob_k'\dromsol_k'(\,\cdot\,,\,\hat\dparamvar_k),\,\cdot\,,\,\hat\dparamvar_k)
  \right]\right|
 \right\}\right).
\end{equation}
The integrand in each term is precisely the integrand of $\Ecal_{2,k}'$ at the
two parameters instances of interest: the trust-region center $\dparamvar_k$
and the candidate step $\hat\dparamvar_k$. Therefore
this refinement process can be repeated iteratively until the conditions on
$\Ecal_{2,k}'$ in (\ref{eqn:inexobj-cond-indiv}) are satisfied. 

\textbf{Reduced-basis construction}
Following
the combined dimension-adaptive greedy method introduced in Section \ref{eq:refinementGradCond} for the gradient
condition, the algorithm interleaves sparse-grid refinement steps with greedy
reduced-basis construction. For a fixed $(\Ical_k',\,\rrob_k')$,
define $\dstochvar^* \in \dstochsp_{\Ical_k'\cup\Ncal(\Ical_k')}$ and
$\dparamvar^* \in \{\dparamvar_k,\,\hat\dparamvar_k\}$
as the quantities that maximize the weighted residual-based error indicator,
i.e.,
\begin{equation} \label{eqn:argmax-greedy2}
	(\dstochvar^*,\,\dparamvar^*) =
 \underset{\substack{\dstochvar\in\dstochsp_{\Ical_k'\cup\Ncal(\Ical_k')}, \\
                     \dparamvar\in\{\dparamvar_k,\,\hat\dparamvar_k\},}}
          {\arg\max} \probdens(\dstochvar)
 \norm{\dpderes(\rrob_k'\dromsol_k'(\dstochvar,\,\dparamvar),\,
                \dstochvar,\,\dparamvar)}_\Thetabold.
\end{equation}
The reduced basis $\rrob_k'$ is updated according to
$\rrob_k' \leftarrow
 \begin{bmatrix} \rrob_k' &
                 \dpdesolSol(\dstochvar^*,\,\dparamvar^*) &
                 \dadjvarSol(\dstochvar^*,\,\dparamvar^*)
 \end{bmatrix}$;
an optional orthogonalization step is usually used to ensure the
reduced basis is full rank and the resulting reduced-order model
is well-conditioned. The argument of the maximization problem in
(\ref{eqn:argmax-greedy2}) is precisely the integrand of $\Ecal_{1,k}'$.
Assuming a minimum-residual reduced-order model is used,
the terms $\Ecal_{1,k}'(\dparamvar_k)$ and $\Ecal_{1,k}'(\hat\dparamvar_k)$
will monotonically decrease with each greedy iteration; the greedy algorithm
proceeds until the requirements on $\Ecal_{1,k}'$ in inequalities (\ref{eqn:inexobj-cond-indiv})
are satisfied. 

\textbf{Refinement algorithm}
The proposed algorithm alternates between sparse-grid refinement
and reduced-basis enrichment exactly as in
Algorithm~\ref{alg:refine-grad}: for a fixed sparse grid, the greedy
method is applied to enrich the reduced basis, then the
reduced-order model is fixed and the sparse grid is refined. The combined
algorithm terminates when all conditions in (\ref{eqn:inexobj-cond-indiv}) are
satisfied.

Algorithm~\ref{alg:refine-obj} summarizes the proposed algorithm.
Similarly to Algorithm~\ref{alg:refine-grad}, this algorithm \emph{refines} a
given sparse-grid/reduced-basis pair and thus implicitly requires
initialization of each quantity. At any iteration $k$, the pair
$(\Ical_k,\,\rrob_k)$ generated by Algorithm~\ref{alg:refine-grad} 
to satisfy the gradient-condition inequalities (\ref{eqn:inexobj-cond-indiv})
at $\dparamvar_k$ is used to initialize Algorithm~\ref{alg:refine-obj}.
If the gradient-condition inequalities turn out to be more restrictive
than that in the objective-condition inequalities (\ref{eqn:inexobj-cond0}),
the algorithm will not modify the sparse grid or reduced basis, i.e.,
$\Ical_k' = \Ical_k$ and $\rrob_k' = \rrob_k$. In this case,
the actual-to-predicted ratio is unity and acceptance of the step is
guaranteed.
%\DPK{Is there any benefit to having separate sparse grid, ROM pairs for
%the gradient and objective function?}
%\MJZ{The refinement is based on different error criteria so, even though
%     I set them to be the same initially, i.e., $\Phibold_k$ is initialized
%     as $\Phibold_{k-1}'$, the refinement modifies the sparse grid, ROM
%     pair for the gradient, objective differently.}

\begin{algorithm}
 \caption{Refine reduced basis and sparse grid for objective condition}
 \label{alg:refine-obj}
 \begin{algorithmic}[1]
 \REQUIRE $\Ical_k$, $\rrob_k$, $\dparamvar_k$, $\hat\dparamvar_k$,
					$\{r_k\}_{k=1}^\infty$ such that $r_k \rightarrow 0$,
          $\omega \in (0,\,1)$, $\alpha_1 > 0$, $\alpha_2 > 0$,
          $\eta < \min\{\eta_1,1-\eta_2\}$
 \ENSURE $\Ical_k'$, $\rrob_k'$
 \STATE \textbf{Set}: $\Ical_k' \leftarrow \Ical_k$,
                      $\rrob_k' \leftarrow \rrob_k$
 \WHILE{
  \begin{equation*}
   \Ecal_{2,k}(\dparamvar_k) + \Ecal_{2,k}(\hat\dparamvar_k) >
   \frac{1}{2\alpha_2}\left(\eta\min\left\{
  \Ebb_{\Ical_k'}\left[\qoi(\rrob_k'\dromsol_k'(\,\cdot\,,\,\dparamvar_k),\,
                            \cdot\,,\,\dparamvar_k)\right] -
  \Ebb_{\Ical_k'}\left[\qoi(\rrob_k'\dromsol_k'(\,\cdot\,,\,\hat\dparamvar_k),\,
                            \cdot\,,\,\hat\dparamvar_k)\right]
  ,\,r_k\right\}\right)^{1/\omega}
 \end{equation*} }
 \vspace{1.5mm}
 \STATE \textbf{Refine index set}: Add index set with largest contribution to
         truncation error
\begin{equation*}
 \Ical_k' \leftarrow \Ical_k' \cup \{\ibm^*\} \quad\text{where}\quad
 \ibm^* = \underset{\ibm \in \Ncal(\Ical_k')}{\arg\max}
     \left(\max\left\{
     \left|\Delta^\ibm \qoi(\rrob_k'\dromsol_k'(\,\cdot\,,\,\dparamvar_k),\,
                            \cdot\,,\,\dparamvar_k)\right|,\,
     \left|\Delta^\ibm \qoi(\rrob_k'\dromsol_k'(\,\cdot\,,\,\hat\dparamvar_k),\,
                            \cdot\,,\,\hat\dparamvar_k)\right|
     \right\}\right)
\end{equation*}
 \WHILE{
  \begin{equation*}
   \Ecal_{1,k}(\dparamvar_k) + \Ecal_{1,k}(\hat\dparamvar_k) >
   \frac{1}{2\alpha_1}\left(\eta\min\left\{
  \Ebb_{\Ical_k'}\left[\qoi(\rrob_k'\dromsol_k'(\,\cdot\,,\,\dparamvar_k),\,
                            \cdot\,,\,\dparamvar_k)\right] -
  \Ebb_{\Ical_k'}\left[\qoi(\rrob_k'\dromsol_k'(\,\cdot\,,\,\hat\dparamvar_k),\,
                            \cdot\,,\,\hat\dparamvar_k)\right]
  ,\,r_k\right\}\right)^{1/\omega}
 \end{equation*} }
 \vspace{1.5mm}
 \STATE \textbf{Evaluate error indicator}: Greedily select
        $\dstochvar \in \dstochsp_{\Ical_k'\cup\Ncal(\Ical_k')}$,
        $\dparamvar \in \{\dparamvar_k,\,\hat\dparamvar_k\}$
        with the largest error
 \begin{equation*}
  \dstochvar^*,\,\dparamvar^* =
  \underset{\substack{\dstochvar \in \dstochsp_{\Ical_k'\cup\Ncal(\Ical_k')}
                      \vspace{0.5mm}\\
                      \dparamvar \in \{\dparamvar_k,\,\hat\dparamvar_k\}}}
           {\arg\max}
           \probdens(\dstochvar)
           \norm{\dpderes(\rrob_k'\dromsol_k'(\,\cdot\,,\,\dparamvar),\,
                          \cdot\,,\,\dparamvar)}
 \end{equation*}
 \STATE \textbf{Reduced-order model construction}: Update reduced basis with
        new snapshot
 \begin{equation*}
  \rrob_k' \leftarrow
  \begin{bmatrix}
   \rrob_k' & \dpdesolSol(\dstochvar^*,\,\dparamvar^*) &
              \dadjvarSol(\dstochvar^*,\,\dparamvar^*)
  \end{bmatrix}
 \end{equation*}
 \ENDWHILE
\ENDWHILE
\end{algorithmic}
\end{algorithm}

\section{Numerical results: optimal boundary control of the incompressible
         Navier--Stokes equations}
\label{sec:num-exp}
In this section, we apply the proposed 
method based on sparse grids and reduced-order models to solve an
optimal control problem governed by the steady-state incompressible
Navier--Stokes equations; we also compare the method's computational efficiency to
that of
existing methods. 

Let $\Omega_x \subset \Rbb^2$ be the channel with a backward facing
step shown in Figure~\ref{fig:bckstp_geom}. The goal of the control problem is
to minimize the expected vorticity in the region immediately downstream of the
step, denoted by $\Omega^*$, by controlling the inflow velocity along the
boundary $\Gamma_c$ subject
to uncertainty both in the fluid viscosity $\nu$ and in the inlet velocity along
boundary $\Gamma_i$. The optimization problem takes the form
\begin{equation}\label{eqn:ins-stoch-opt}
\begin{aligned}
& \underset{\dparamvar \in \Rbb^\nmu}
           {\mathrm{minimize}}
& & \int_{\dstochsp} \probdens(\dstochvar)
    \left[\frac{1}{2}\int_{\Omega^*}
                        |\nabla \times u(x,\,\dstochvar,\,\mubold)|^2\,d\Omega +
  \frac{\alpha}{2}\int_{\Gamma_c}
                        |g(x,\,\dparamvar)|^2\,d\Gamma\right]d\dstochvar,
\end{aligned}
\end{equation}
where $u(x,\,\dstochvar,\,\dparamvar)$ is the solution of the incompressible
Navier--Stokes equations
\begin{equation} \label{eqn:ins}
 \begin{aligned}
  -\nu\Delta u + (u\cdot\nabla u)u + \nabla p = 0 ,&\quad x \in \Omega \\
  \nabla \cdot u = 0 ,&\quad x \in \Omega \\
  (\nabla u - pI) n = 0 ,&\quad x \in \Gamma_o \\
  u = g ,&\quad x \in \Gamma_c \\
  u = b ,&\quad x \in \Gamma_i \\
  u = 0 ,&\quad x \in \partial\Omega \setminus
                     (\Gamma_i \cup \Gamma_c \cup \Gamma_o).
 \end{aligned}
\end{equation}
\ifbool{fastcompile}{}{
\begin{figure}
 \centering
 \begin{tikzpicture}
\begin{axis}[
axis equal image,
axis line style={gray},
axis x line*=bottom,
axis y line*=left,
ytick = {0.0, 0.5, 1.0},
xmin=0.0,
width=0.8\textwidth,
ymin=0.0,
ymax=1.0,
xmax=8.0]
\addplot [fill opacity=0.5, fill=black!30!white, forget plot]
coordinates {
(  0.00000000,   0.50000000)
(  1.00000000,   0.50000000)
(  1.00000000,   0.00000000)
(  8.00000000,   0.00000000)
(  8.00000000,   1.00000000)
(  0.00000000,   1.00000000)
(  0.00000000,   0.50000000)};

\addplot [thin, red!40!white!, dashed, fill opacity=0.3, fill=red!40!white, forget plot]
coordinates {
(  1.00000000,   0.50000000)
(  1.00000000,   0.00000000)
(  3.00000000,   0.00000000)
(  3.00000000,   0.50000000)
(  1.00000000,   0.50000000)};

\addplot [thick, color=blue]
coordinates {
(  0.00000000,   1.00000000)
(  0.00000000,   0.50000000)};\label{line:bckstp0:inflow_stoch}

\addplot [thick, color=black]
coordinates {
(  0.00000000,   0.50000000)
(  1.00000000,   0.50000000)};\label{line:bckstp0:wall}

\addplot [thick, color=red]
coordinates {
(  1.00000000,   0.50000000)
(  1.00000000,   0.00000000)};\label{line:bckstp0:inflow_param}

\addplot [thick, color=black, forget plot]
coordinates {
(  1.00000000,   0.00000000)
(  8.00000000,   0.00000000)};

\addplot [thick, color=yellow]
coordinates {
(  8.00000000,   0.00000000)
(  8.00000000,   1.00000000)};\label{line:bckstp0:outflow}

\addplot [thick, color=black, forget plot]
coordinates {
(  8.00000000,   1.00000000)
(  0.00000000,   1.00000000)};

\node[left]    at    (axis cs:8.0, 0.5) {$\Gamma_o$};
\node[left]    at    (axis cs:1.0, 0.25) {$\Gamma_c$};
\node[right]    at    (axis cs:0.0, 0.75) {$\Gamma_i$};
\node[]    at    (axis cs:2.0, 0.25) {$\Omega^*$};
\end{axis}
\end{tikzpicture}
 \caption{Geometry and boundary conditions for backward facing step.
          Boundary conditions: viscous wall (\ref{line:bckstp0:wall}),
          parametrized inflow (\ref{line:bckstp0:inflow_param}),
          stochastic inflow (\ref{line:bckstp0:inflow_stoch}),
          outflow (\ref{line:bckstp0:outflow}). Vorticity magnitude is minimized
          in the red shaded region.}
 \label{fig:bckstp_geom}
\end{figure}
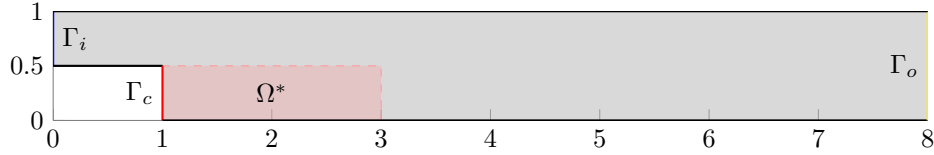
}
The fluid viscosity and inflow are uncertain and parametrized with a
two-dimensional stochastic space $\dstochsp = [-1,\,1]^2$, as
\begin{equation*}
 \nu(\dstochvar)^{-1} = 50(1-y_1) + 250(y_1+1), \qquad
 b(x,\,\dstochvar) = (8+y_2)(x_2-0.5)(1-x_2)
\end{equation*}
with a uniform probability distribution
$\probdens(\dstochvar)d\dstochvar = 2^{-2}d\dstochvar$.
In all numerical experiments, the partial differential equation in
(\ref{eqn:ins}) is discretized with $232$ $\Qbb_2-\Qbb_1$ Taylor-Hood
elements, yielding a state space of dimension $\nU = 2034$, after application of
the essential boundary conditions. The velocity at all nodes along the
parametrized inlet $\Gamma_c$, are taken as optimization variables,
yielding a parameter space of dimension $\nmu = 38$, and the finite element
basis functions restricted to $\Gamma_c$ define the boundary control
$g(x,\,\dparamvar)$. The control regularization parameter is set to
$\alpha = 0.1$.
The PDE in this example is quadratic in $u$ and linear in $\mubold$.
Similarly, the objective is quadratic in both $u$ and $\mubold$. Therefore, the partial
derivatives of the objective function and PDE are at most linear in $u$ and $\mubold$,
which establishes the assumptions required in Appendix~\ref{app:resbnd} (bounded, linear
mappings are Lipschitz continuous).

The flow velocity at the parametrized inlet is initially set to zero
($\dparamvar = 0$), i.e., a viscous wall, which leads to a recirculation
in the mean flow and standard-deviation offset defined as
\begin{equation*}
	\bar{u}(x,\,\dparamvar) = \Ebb[u(x,\,\cdot\,,\,\dparamvar)]
	\qquad\text{and}\qquad
\bar{u}_\pm(x,\,\dparamvar) = \bar{u}(x,\,\dparamvar) \pm
            \sqrt{\Ebb[(u(x,\,\cdot\,,\,\dparamvar)-\bar{u}(x,\,\dparamvar))^2]},
\end{equation*}
respectively (see Figure~\ref{fig:ins:hdmopt:stochiso:vort-quiver}). The local
solution of the stochastic optimal control problem shown in
Figure~\ref{fig:ins:hdmopt:stochiso:vort-quiver} effectively 
eliminates the recirculation region in not only the mean flow, but also in the
standard-deviation offsets.
\ifbool{fastcompile}{}{
\begin{figure}
 \centering
 \input{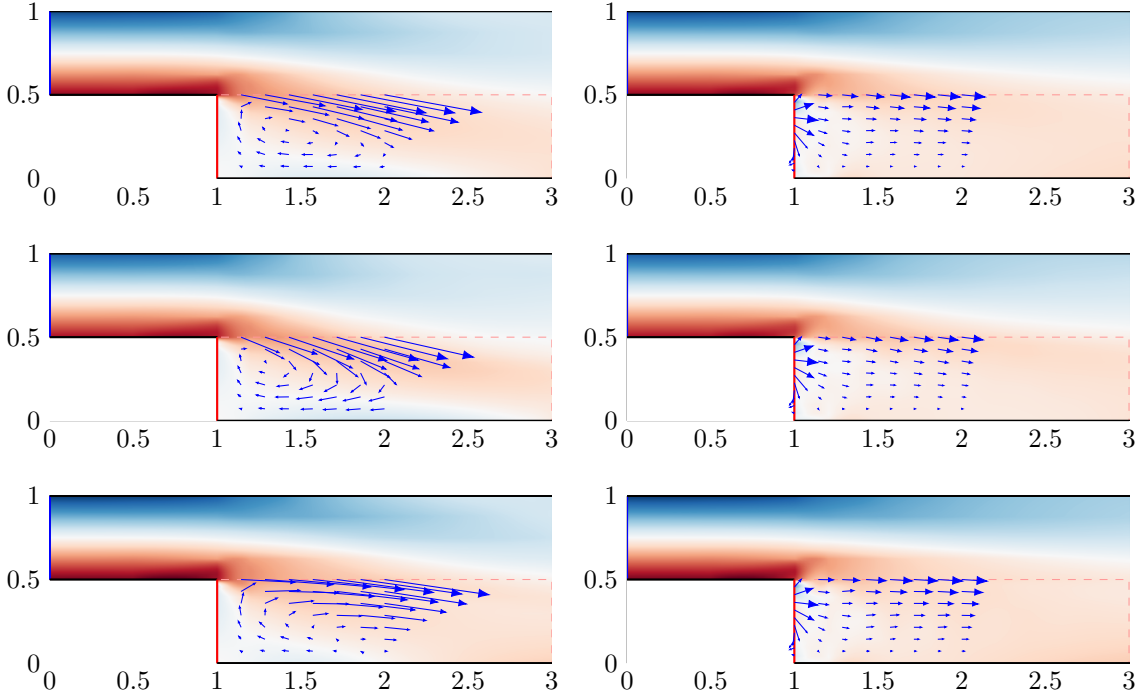}
 \caption{The mean flow $\bar{u}(x,\,\dparamvar)$ (\emph{top}) and
          standard deviation offsets $\bar{u}_{-}(x,\,\dparamvar)$
          (\emph{center}), $\bar{u}_{+}(x,\,\dparamvar)$ (\emph{bottom})
          corresponding to the uncontrolled, $\dparamvar = \zerobold$,
          (\emph{left}) and controlled flow (\emph{right}). Boundary
          control along $\Gamma_c$ effectively eliminates the recirculation
          region.}
 \label{fig:ins:hdmopt:stochiso:vort-quiver}
\end{figure}
}

We apply the proposed method  to solve problem
(\ref{eqn:ins-stoch-opt}). We construct sparse grids using
one-dimensional, nested Clenshaw--Curtis quadrature rules
\cite{clenshaw1960method}, which ensure HDM and ROM evaluations can be
recycled between sparse-grid refinement levels because the quadrature nodes at
level $j$ are a subset of those at level $j+1$.
We employ minimum-residual primal and adjoint reduced-order models because the
PDE operators are not symmetric positive definite. For simplicity and
implementation efficiency, the primal and adjoint reduced-order models
are constructed to minimize their respective residuals in the Euclidean norm,
i.e., $\Thetabold = \Thetabold^\dadjvar = \Ibm$. Because we use minimum-residual
primal and adjoint reduced-order models,
Algorithms~\ref{alg:refine-grad}--\ref{alg:refine-obj} are guaranteed to
return sparse grids and reduced bases that satisfy
the gradient condition
(\ref{eqn:inexgrad-cond0}) and objective condition (\ref{eqn:inexobj-cond0}), and global convergence
is thus ensured, % (Table~\ref{tab:ins-trconv-mii}), 
provided the regularity conditions
in Appendix~\ref{app:resbnd} are satisfied. The reduced basis is initialized
using the primal and adjoint solution and primal solution sensitivity at
$\dstochvar=\dstochvar_0 = \zerobold$ and $\dparamvar=\dparamvar_0$, i.e.,
\begin{equation*}
\rrob_{-1}' = \begin{bmatrix}
              \dpdesolSol(\dstochvar_0,\,\dparamvar_0) &
              \displaystyle{
                  \pder{\dpdesolSol}{\dparamvar}(\dstochvar_0,\,\dparamvar_0)} &
              \dadjvarSol(\dstochvar_0,\,\dparamvar_0)
             \end{bmatrix}
\end{equation*}
to enhance the parametric robustness of the reduced-order model at early
trust-region iterations with few HDM samples. Subsequently, the basis
is updated using only primal and adjoint solutions as prescribed by 
Algorithms~\ref{alg:refine-grad} and \ref{alg:refine-obj}, because
including $\nmu = 38$ sensitivity solutions at every HDM sample
would result in a computationally inefficient method.
An approximate solution of the trust
region subproblem is computed using the Steihaug--Toint Conjugate Gradient
method \cite{toint1981towards, steihaug1983conjugate, conn2000trust}
to minimize a quadratic approximation of $\aprxmod_k(\dparamvar)$ in (\ref{eq:approximationModels}) within the trust region. 

Given the numerous user-defined constants used in the proposed method that can
impact the performance of the scheme (convergence is guaranteed regardless of their
values, provided they satisfy the conditions discussed previously), we discuss
general strategies for selecting these parameters and the specific choices made for the
problem at hand. The parameter $\kappa_\varphi$ gives the user indirect control of the refinement of the sparse grid and reduced basis at the initial iteration, and therefore all subsequent iterations. Small values of $\kappa_\varphi$ would require a refined sparse grid and reduced basis at the initial guess, while larger values would accept a coarser sparse grid and reduced basis. In this work, we simply take $\kappa_\varphi = 1$ so the refinement of the approximation model is solely determined from the scaling between the gradient error indicator and model gradient.
The parameters $\eta_1, \eta_2$ are standard trust-region parameters that
determine when to accept a particular candidate step and modify the trust-region
radius based on the value of the actual-to-predicted reduction, and $\gamma$ is
the factor used to modify the trust-region
radius. In this work, we make the standard choices $\eta_1 = 0.1$, $\eta_2 = 0.75$, and
$\gamma = 0.5$ \cite{kouri2013trust}. The parameters $\eta, \omega$, and the forcing
sequence $\{r_k\}_{k=0}^\infty$ were introduced in \cite{kouri2014inexact} to define the error condition on the objective required for convergence. We use the values $\eta = 0.1$, $\omega = 0.1$, and $r_k = 1/(k+1)$ suggested in that paper.
In practice, we take $\beta_1, \beta_3, \beta_4$ such that the three terms in the
error indicator (\ref{eqn:gradbnd0}) are approximately equal when the algorithm is
initialized, i.e.,
$\beta_1 \approx \beta\Ecal_{1,k}(\mubold_0)^{-1}$,
$\beta_3 \approx \beta\Ecal_{3,k}(\mubold_0)^{-1}$,
$\beta_4 \approx \beta\Ecal_{4,k}(\mubold_0)^{-1}$, where $k = -1$ and
$\beta > 0$ gives indirect control over the refinement of the approximation model
at the initial guess. In this work, we take $\beta_1 = \beta_3 = \beta_4 = 1$ because
the error terms have similar scale and we use the natural scaling between the
error indicator terms and the model gradient to determine the refinement of the
initial approximation model.
Similarly, we take $\alpha_1, \alpha_2$ such that the two terms in the
error indicator (\ref{eqn:objbnd0}) are approximately equal when the algorithm is initialized, i.e.,
$\alpha_1 \approx \alpha\Ecal_{1,k}'(\mubold_0)^{-1}$ and
$\alpha_2 \approx \alpha\Ecal_{2,k}'(\mubold_0)^{-1}$,
where $k = -1$ and $\alpha > 0$ gives indirect control over step acceptance.
In this work, we take $\alpha_1 = \alpha_2 = 10^{-2}$ because
the error terms have similar scale and we aim to promote step acceptance in early
iterations.

%The trust region
%algorithm parameters, e.g., those used in Algorithm~\ref{alg:tr-iog-gen},
%are taken as: $\eta_1 = 0.1$, $\eta_2 = 0.75$, $\eta = 0.1$, $\gamma = 0.5$,
%$r_k = 1/(k+1)$, and $\omega = 0.9$. An approximate solution of the trust
%region subproblem is computed using the Steihaug--Toint Conjugate Gradient
%method \cite{toint1981towards, steihaug1983conjugate, conn2000trust}
%to minimize a quadratic approximation of $\aprxmod_k(\dparamvar)$ in
%(\ref{eq:approximationModels}) within the trust region. The constants used
%to weight the contributions of the individual terms in the gradient error
%indicator are taken as unity, $\beta_1 = \beta_3 = \beta_4 = 1.0$, and
%$\kappa_\varphi = 1.0$. To equally weight the terms in the objective
%function error indicator and weaken the error condition on the objective
%function, i.e., to make step acceptance more likely, the constants are
%taken as $\alpha_1 = \alpha_2 = 0.01$. 

%\DPK{Again, convergence is only guaranteed if the solution is sufficiently
%regular with respect to the uncertain parameters.  I don't think we have
%verified this for Navier--Stokes...maybe someone has???}
%\MJZ{I softened it with ``provided the regularity conditions...''}.

\begin{figure}
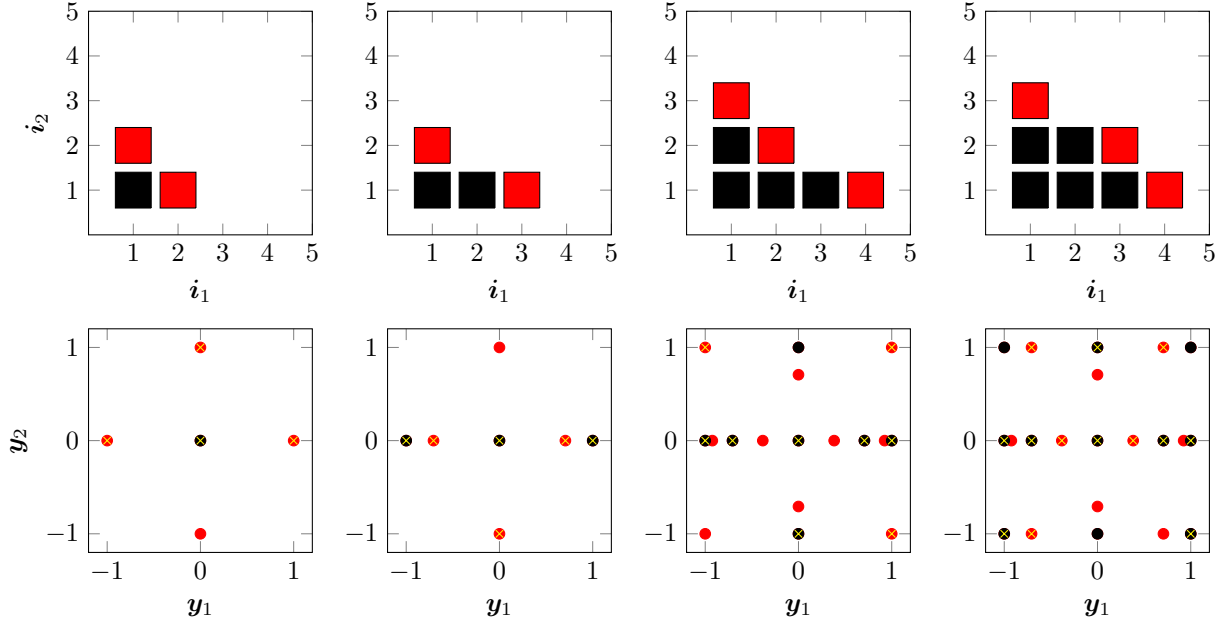

 \centering
 \begin{tikzpicture}
\begin{groupplot}[
    group style={
        group size=4 by 2,
        horizontal sep=1.0cm,
        vertical sep=1.25cm
    },
    axis equal image,
    width=0.32\textwidth
]

\nextgroupplot[xlabel={$\ibm_1$}, ylabel={$\ibm_2$},
               xmin=0, xmax=5, ymin=0, ymax=5,
               xtick={1, 2, 3, 4, 5}, ytick={1, 2, 3, 4, 5}]
\input{tikz/ins/romopt_aniso_iog_prim-sens-dual_ctr-tcg.sg.iter00.idxset.tikz}

\nextgroupplot[xlabel={$\ibm_1$},
               xmin=0, xmax=5, ymin=0, ymax=5,
               xtick={1, 2, 3, 4, 5}, ytick={1, 2, 3, 4, 5}]
\input{tikz/ins/romopt_aniso_iog_prim-sens-dual_ctr-tcg.sg.iter03.idxset.tikz}

\nextgroupplot[xlabel={$\ibm_1$},
               xmin=0, xmax=5, ymin=0, ymax=5,
               xtick={1, 2, 3, 4, 5}, ytick={1, 2, 3, 4, 5}]
\input{tikz/ins/romopt_aniso_iog_prim-sens-dual_ctr-tcg.sg.iter05.idxset.tikz}

\nextgroupplot[xlabel={$\ibm_1$},
               xmin=0, xmax=5, ymin=0, ymax=5,
               xtick={1, 2, 3, 4, 5}, ytick={1, 2, 3, 4, 5}]
\input{tikz/ins/romopt_aniso_iog_prim-sens-dual_ctr-tcg.sg.iter06.idxset.tikz}

\nextgroupplot[xlabel={$\ybm_1$}, ylabel={$\ybm_2$}]
\input{tikz/ins/romopt_aniso_iog_prim-sens-dual_ctr-tcg.sg.iter00.nodeset.tikz}

\nextgroupplot[xlabel={$\ybm_1$}]
\input{tikz/ins/romopt_aniso_iog_prim-sens-dual_ctr-tcg.sg.iter03.nodeset.tikz}

\nextgroupplot[xlabel={$\ybm_1$}]
\input{tikz/ins/romopt_aniso_iog_prim-sens-dual_ctr-tcg.sg.iter05.nodeset.tikz}

\nextgroupplot[xlabel={$\ybm_1$}]
\input{tikz/ins/romopt_aniso_iog_prim-sens-dual_ctr-tcg.sg.iter06.nodeset.tikz}

\end{groupplot}
\end{tikzpicture}
 \caption{The sparse-grid index set (\emph{top}) and corresponding quadrature
          nodes (\emph{bottom}) used by the proposed method at trust-region
          iterations $k=0,\,3,\,5,\,6$ (\emph{left-to-right}).
          Legend: sparse grid $\Ical_k$ (\ref{line:sg:idx}),
                  sparse grid neighbors $\Ncal(\Ical_k)$(\ref{line:sg:neigh}).
          The points in stochastic space where the primal and adjoint HDM
          are sampled to construct the reduced basis $\rrob_k$ are indicated
					with a (\ref{line:sg:rom}).}
 \label{fig:ins:sg}
\end{figure}

Table~\ref{tab:ins-trconv-mii} reports the convergence history of the
proposed method. The method converges to a first-order critical
point ($\norm{\nabla \obj(\dparamvar_k)} \rightarrow 0$); after only $8$
trust-region iterations the first-order optimality condition has reduced nearly
$4$ orders of magnitude from the initial (sub-optimal) control. At early
iterations, the approximation model $\aprxmod_k(\dparamvar)$ and true
objective $\obj(\dparamvar)$ agree only to one digit at the candidate
step $\hat\dparamvar_k$; this is expected because a coarse sparse grid
and only a few additional primal and adjoint HDM samples are required
(Figure~\ref{fig:ins:sg}) to meet the gradient condition. However,
despite the inexpensive approximation with limited accuracy, the resulting
approximation model is sufficiently accurate
to significantly reduce the \emph{true} objective function during the
trust-region iteration. As the algorithm approaches a local minimum, the 
the gradient condition (\ref{eqn:inexgrad-cond0}) places more stringent requirements on the model
error and, as a result, the approximation model consists of a finer sparse
grid with more HDM samples used to construct the reduced basis
(Figure~\ref{fig:ins:sg}) and thus provides a better approximation of the objective
function.
\ifbool{fastcompile}{}{
\begin{table}[htb]
 \caption{Convergence history of the proposed method applied to the optimal
          control of the incompressible Navier--Stokes equation in
          (\ref{eqn:ins-stoch-opt}).}
 \label{tab:ins-trconv-mii}
 \begin{tabular}{cccccccc}
  \toprule
  $\obj(\dparamvar_k)$ & $\aprxmod_k(\dparamvar_k)$ &
  $\obj(\hat\dparamvar_k)$ & $\aprxmod_k(\hat\dparamvar_k)$ &
  $\norm{\nabla \obj(\dparamvar_k)}$ & $\varrho_k$ &
  $\trrad_k$ & Success? \\\midrule
  1.0740e+00  &  1.0805e+00  &  8.4412e-01  &  8.6172e-01  &  1.8723e+00  &  1.0000e+00  &  2.500e-01  &  True \\
8.4412e-01  &  8.4351e-01  &  7.4896e-01  &  7.4628e-01  &  1.3292e+00  &  1.0000e+00  &  5.000e-01  &  True \\
7.4896e-01  &  7.3757e-01  &  7.3766e-01  &  7.2654e-01  &  3.3224e-01  &  8.6570e-01  &  1.000e+00  &  True \\
7.3766e-01  &  7.3429e-01  &  7.3601e-01  &  7.3204e-01  &  1.1425e-01  &  7.3229e-01  &  2.000e+00  &  True \\
7.3601e-01  &  7.3250e-01  &  7.3548e-01  &  7.3207e-01  &  7.9688e-02  &  1.2288e+00  &  2.000e+00  &  True \\
7.3548e-01  &  7.3625e-01  &  7.3536e-01  &  7.3612e-01  &  1.4001e-02  &  9.4543e-01  &  4.000e+00  &  True \\
7.3536e-01  &  7.3549e-01  &  7.3535e-01  &  7.3549e-01  &  4.5723e-03  &  1.1648e+00  &  8.000e+00  &  True \\
7.3535e-01  &  7.3549e-01  &         -  &         -  &  7.8794e-04  &         -  &         -  &         - \\

  \bottomrule
 \end{tabular}
\end{table}
}

We assess the computational performance of the proposed method using two baseline methods for comparison. The first method uses a fixed, $5$-level isotropic
sparse grid to perform integration in stochastic space using only HDM
evaluations and the optimization problem is solved with the BFGS algorithm
\cite{nocedal2006numerical}. The second method uses
the trust-region dimension-adaptive sparse grid approach proposed in
Ref.~\cite{kouri2014inexact}.
For brevity, we abbreviate these methods SG-ISO and SG-TR, respectively.
We refer to our proposed trust-region approach with dimension-adaptive
sparse grids and reduced-order models as SG-ROM-TR.

The SG-TR method requires a large number of primal and adjoint HDM
evaluations as the trust-region iterations progress, because the evaluation
of the models and their error indicators rely solely on HDM evaluations,
albeit on a highly adapted sparse grid. In contrast, our SG-ROM-TR method
requires HDM evaluations only to construct the reduced basis; evaluations of
the trust-region models and their error indicators rely only on reduced-order
model evaluations. As a result, the proposed SG-ROM-TR method requires one to
two orders of magnitude fewer queries to the HDM
(Figure~\ref{fig:ins:hdm-query-vs-majit}) at the cost of a large number of
ROM queries for various reduced-basis dimensions, $\kU \in [39,\,165]$
(Figure~\ref{fig:ins:rom-query-vs-majit}).
\ifbool{fastcompile}{}{
\begin{figure}
 \centering
 \begin{tikzpicture}
\begin{groupplot}[
  group style={
    group size=2 by 1, 
    horizontal sep=2.4cm
  },
  width=7.5cm,
  height=5.5cm,
  xmin=0, xmax=8,
  ymin=1e0, ymax=1e4,
  ymode=log]

\nextgroupplot[xlabel={Major iterations}, ylabel={Primal HDM evaluations}]
\addplot [blue, dashed, thick, mark=square*, mark size=2, mark options={solid}]  table[x index=0, y index=15, select coords between index={0}{6}] {dat/ins/hdmopt_stoch_aniso_kouri14.dat}; \label{line:ins:kouri14}
\addplot [black, solid, thick, mark=*, mark size=2, mark options={solid}]  table[x index=0, y index=15] {dat/ins/romopt_aniso_iog_prim-sens-dual_ctr-tcg.dat}; \label{line:ins:romopt-iog}

\nextgroupplot[xlabel={Major iterations}, ylabel={Adjoint HDM evaluations}]
\addplot [blue, dashed, thick, mark=square*, mark size=2, mark options={solid}]  table[x index=0, y index=16, select coords between index={0}{6}] {dat/ins/hdmopt_stoch_aniso_kouri14.dat};
\addplot [black, solid, thick, mark=*, mark size=2, mark options={solid}]  table[x index=0, y index=16] {dat/ins/romopt_aniso_iog_prim-sens-dual_ctr-tcg.dat};

\end{groupplot}
\end{tikzpicture}
 \caption{Cumulative number of HDM primal and adjoint evaluations as
          the trust-region iterations in the various trust region algorithms
					progress: dimension-adaptive sparse grid (SG-TR)
          \cite{kouri2014inexact} (\ref{line:ins:kouri14}) and proposed
          method (ROM-SG-TR) (\ref{line:ins:romopt-iog}).}
 \label{fig:ins:hdm-query-vs-majit}
\end{figure}
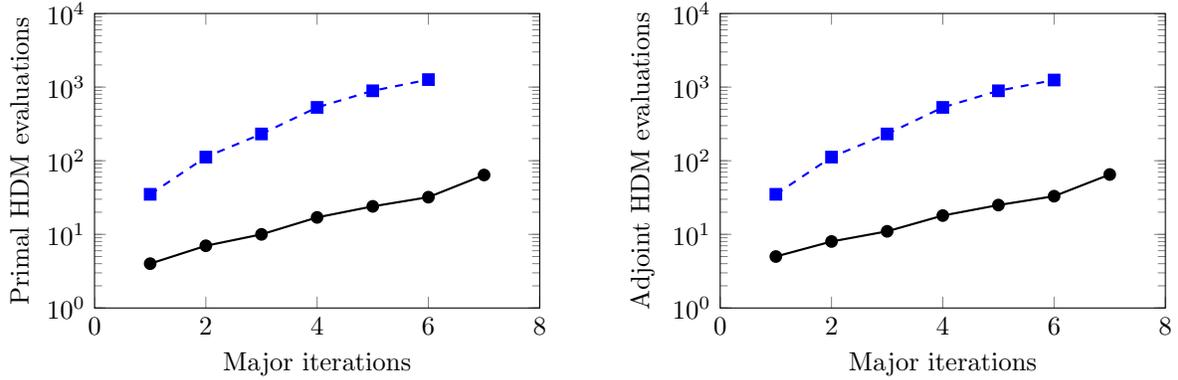
}
\ifbool{fastcompile}{}{
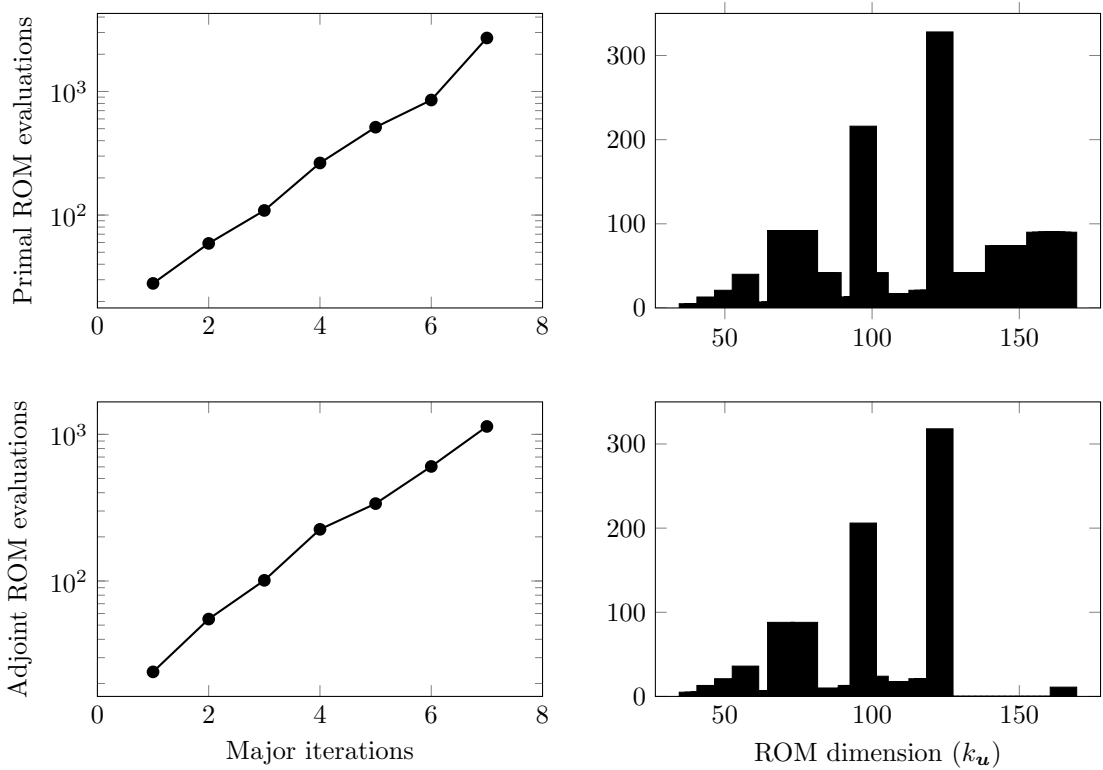
\begin{figure}
 \centering
 \begin{tikzpicture}
\begin{groupplot}[
    group style={
        group size=2 by 2,
        horizontal sep=1.5cm,
        vertical sep=1.25cm
    },
    width=7.5cm,
    height=5.5cm,
]

\nextgroupplot[xmin=0, xmax=8, ymode=log,
               ylabel={Primal ROM evaluations}]
\addplot [black, solid, thick, mark=*, mark size=2, mark options={solid}]  table[x index=0, y index=17] {dat/ins/romopt_aniso_iog_prim-sens-dual_ctr-tcg.dat};

\nextgroupplot[ybar, fill, ymin=0, ymax=350]
\addplot[ybar, fill]  table[x index=0, y index=1] {dat/ins/romopt_aniso_iog_prim-sens-dual_ctr-tcg.glob.dat};

\nextgroupplot[xmin=0, xmax=8, ymode=log,
               xlabel={Major iterations},
               ylabel={Adjoint ROM evaluations}]
\addplot [black, solid, thick, mark=*, mark size=2, mark options={solid}]  table[x index=0, y index=18] {dat/ins/romopt_aniso_iog_prim-sens-dual_ctr-tcg.dat};

\nextgroupplot[ymin=0, ymax=350, xlabel={ROM dimension ($\kU$)}]
\addplot[ybar, fill]  table[x index=0, y index=3] {dat/ins/romopt_aniso_iog_prim-sens-dual_ctr-tcg.glob.dat};
\end{groupplot}
\end{tikzpicture}
 \caption{Number of primal and adjoint reduced-order model evaluations
          required by proposed method. \emph{Left}: Cumulative number
          evaluations as the trust-region iterations progress. \emph{Right}:
					Histogram of the number of evaluations as a function of 
          reduced-basis dimension
          $k_\ubm$.}
 \label{fig:ins:rom-query-vs-majit}
\end{figure}
}

Because the proposed and baseline methods have different costs (i.e.,
the proposed method requires only HDM and ROM evaluations, while the baseline
methods employ only HDM evaluations), care must be taken when
assessing their performance. The ultimate cost metric of interest
is wall time; however, it is well-known that hyperreduction
\cite{barrault2004empirical} must be incorporated to reduce the
complexity associated with evaluating the nonlinear terms. Integration of
hyperreduction in the proposed framework is the subject of ongoing work
so we turn to another error metric. To assess the speedups that can be
realized by this method, the following simple cost model is introduced
\begin{equation} \label{eqn:ins:cost}
 C = n_{hp}+n_{ha}\bar{n}_h^{-1} + \tau^{-1}(n_{rp}+n_{ra}\bar{n}_r^{-1})
\end{equation}
where $C$ is the total cost associated with a particular method in the
units of \emph{equivalent number of primal HDM queries}, $n_{hp}$ is
the number of primal HDM queries, $n_{ha}$ is the number of adjoint
HDM queries, $\bar{n}_h$ is the average number of nonlinear iterations
required to solve the primal HDM, $n_{rp}$ is the number of primal ROM
queries, $n_{ra}$ is the number of adjoint ROM queries, $\bar{n}_r$ is
the average number of nonlinear iterations required to solve the primal
ROM, and $\tau$ is the ratio of the cost of a primal HDM query to a primal
ROM query. This cost model assumes a primal HDM (ROM) solve is $\bar{n}_h$
($\bar{n}_r$) times as expensive as an adjoint solve and a primal HDM
solve is $\tau$ times as expensive as a primal ROM solve. Numerical
experiments at the initial and optimal controls suggest
$\bar{n}_h = \bar{n}_r = 5$ and a range of $\tau$ values are considered.

As shown in Figure~\ref{fig:ins:obj-vs-cost}, for costly-to-evaluate
reduced-order models ($\tau = 1$), the proposed SG-ROM-TR method
demonstrates faster convergence, in terms of the cost $C$, than the
brute-force approach SG-ISO and similar convergence to the state-of-the-art
SG-TR method. This is expected because SG-TR and SG-ROM-TR are adaptive
algorithms with refinement at each iteration tailored to the convergence
requirements, while the SG-ISO method relies solely on HDM
evaluations on a fine sparse grid for all iterations. For a modest ROM speedup
($\tau = 10$), the computational cost required for the SG-ROM-TR method to
reach a given value of the objective function is nearly an order of magnitude
lower than that of the state-of-the-art SG-TR method. The hypothetical
situation of a free reduced-order model ($\tau = \infty$) shows the speedup
attainable by our method in comparison to SG-TR is at most $500$.
The $\tau$-scenarios are not realistic in the sense that they assume the
speedup of the ROM is constant across all iterations, which cannot be the
case since the size of the ROM varies by more than a factor of $4$ between
the smallest and largest ROM used (Figure~\ref{fig:ins:rom-query-vs-majit}).
However, the scenarios do provide useful information because they bound the true
cost, i.e., if the smallest and largest ROMs considered have speedup factors
of $\tau_0$ and $\tau_1$, respectively, the true cost-objective curve will lie
between the $\tau_0$ and $\tau_1$ curves. For example, in the modest scenario of
$\tau_0 = 100$ and $\tau_1 = 10$, the true cost-objective curve will lie
between the (\ref{line:ins:romopt-iog-cost-2}) and
(\ref{line:ins:romopt-iog-cost-1}) curves in Figure~\ref{fig:ins:obj-vs-cost}.
This will
lead to a significant speedup of the proposed method over both the SG-ISO and SG-TR methods.
\ifbool{fastcompile}{}{
\begin{figure}
 \def\optimObjIns{7.353521063821857329e-01}
\def\nNewtHdm{5} % number Newton evaluations
\def\costHdmEval{1.0/\nNewtHdm}
\def\nQuadIso{65} % number of quadrature points in isotropic SG

\begin{tikzpicture}
\begin{groupplot}[
    group style={
        group size=2 by 1,
        horizontal sep=2.4cm
    },
    xmin=1,
    xmax=1e4,
    xlabel={Cost},
    xmode=log, ymode=log,
    width=7.5cm,
    height=5.5cm
]

\nextgroupplot[ymin=5e-7, ymax=1.2, ylabel={$|\obj(\dparamvar) - \obj(\dparamvar^*)|$}]
\addplot [red, solid, thick]  table[x expr=\costHdmEval*(\nNewtHdm+1)*\nQuadIso*(1+\coordindex), y expr=abs(\thisrowno{1}-\optimObjIns), select coords between index={3}{54}] {dat/ins/hdmopt_stochiso.hist.dat}; \label{line:ins:stochiso}
\addplot [blue, dashed, thick, mark=square*, mark size=2, mark options={solid}, mark repeat={1}]  table[x expr=\costHdmEval*(\nNewtHdm*\thisrowno{15}+\thisrowno{16}), y expr=abs(\thisrowno{1}-\optimObjIns), select coords between index={0}{6}] {dat/ins/hdmopt_stoch_aniso_kouri14.dat};
\addplot [black, solid, thick, mark=x, mark size=2, mark options={solid}, mark repeat={1}]  table[x expr=\costHdmEval*(\nNewtHdm*\thisrowno{15}+\thisrowno{16})+0*\costHdmEval*(\nNewtHdm*\thisrowno{17}+\thisrowno{18}), y expr=abs(\thisrowno{1}-\optimObjIns)] {dat/ins/romopt_aniso_iog_prim-sens-dual_ctr-tcg.dat}; \label{line:ins:romopt-iog-cost-infty}
\addplot [black, solid, thick, mark=*, mark size=2, mark options={solid}, mark repeat={1}]  table[x expr=\costHdmEval*(\nNewtHdm*\thisrowno{15}+\thisrowno{16})+0.01*\costHdmEval*(\nNewtHdm*\thisrowno{17}+\thisrowno{18}), y expr=abs(\thisrowno{1}-\optimObjIns)] {dat/ins/romopt_aniso_iog_prim-sens-dual_ctr-tcg.dat}; \label{line:ins:romopt-iog-cost-2}
\addplot [black, solid, thick, mark=triangle*, mark size=2, mark options={solid}, mark repeat={1}]  table[x expr=\costHdmEval*(\nNewtHdm*\thisrowno{15}+\thisrowno{16})+0.1*\costHdmEval*(\nNewtHdm*\thisrowno{17}+\thisrowno{18}), y expr=abs(\thisrowno{1}-\optimObjIns)] {dat/ins/romopt_aniso_iog_prim-sens-dual_ctr-tcg.dat}; \label{line:ins:romopt-iog-cost-1}
\addplot [black, solid, thick, mark=diamond*, mark size=2, mark options={solid}, mark repeat={1}]  table[x expr=\costHdmEval*(\nNewtHdm*\thisrowno{15}+\thisrowno{16})+\costHdmEval*(\nNewtHdm*\thisrowno{17}+\thisrowno{18}), y expr=abs(\thisrowno{1}-\optimObjIns)] {dat/ins/romopt_aniso_iog_prim-sens-dual_ctr-tcg.dat}; \label{line:ins:romopt-iog-cost-0}

\nextgroupplot[ymin=1e-4, ymax=10, ylabel={$\norm{\nabla \obj(\dparamvar)}$}]
\addplot [red, solid, thick]  table[x expr=\costHdmEval*(\nNewtHdm+1)*\nQuadIso*(1+\coordindex), y index=2, select coords between index={3}{57}] {dat/ins/hdmopt_stochiso.hist.dat};
\addplot [blue, dashed, thick, mark=square*, mark size=2, mark options={solid}, mark repeat={1}]  table[x expr=\costHdmEval*(\nNewtHdm*\thisrowno{15}+\thisrowno{16}), y index=5, select coords between index={0}{6}] {dat/ins/hdmopt_stoch_aniso_kouri14.dat};
\addplot [black, solid, thick, mark=x, mark size=2, mark options={solid}, mark repeat={1}]  table[x expr=\costHdmEval*(\nNewtHdm*\thisrowno{15}+\thisrowno{16})+0*\costHdmEval*(\nNewtHdm*\thisrowno{17}+\thisrowno{18}), y index=5] {dat/ins/romopt_aniso_iog_prim-sens-dual_ctr-tcg.dat};
\addplot [black, solid, thick, mark=*, mark size=2, mark options={solid}, mark repeat={1}]  table[x expr=\costHdmEval*(\nNewtHdm*\thisrowno{15}+\thisrowno{16})+0.01*\costHdmEval*(\nNewtHdm*\thisrowno{17}+\thisrowno{18}), y index=5] {dat/ins/romopt_aniso_iog_prim-sens-dual_ctr-tcg.dat};
\addplot [black, solid, thick, mark=triangle*, mark size=2, mark options={solid}, mark repeat={1}]  table[x expr=\costHdmEval*(\nNewtHdm*\thisrowno{15}+\thisrowno{16})+0.1*\costHdmEval*(\nNewtHdm*\thisrowno{17}+\thisrowno{18}), y index=5] {dat/ins/romopt_aniso_iog_prim-sens-dual_ctr-tcg.dat};
\addplot [black, solid, thick, mark=diamond*, mark size=2, mark options={solid}, mark repeat={1}]  table[x expr=\costHdmEval*(\nNewtHdm*\thisrowno{15}+\thisrowno{16})+\costHdmEval*(\nNewtHdm*\thisrowno{17}+\thisrowno{18}), y index=5] {dat/ins/romopt_aniso_iog_prim-sens-dual_ctr-tcg.dat};

\end{groupplot}
\end{tikzpicture}
 \caption{Convergence of the objective function (\emph{left}) and gradient
         (\emph{right}) as a function of the cost metric in
         (\ref{eqn:ins:cost}) for the proposed SG-ROM-TR method for several
         values of the speedup factor of the reduced-order model:
         $\tau = 1$ (\ref{line:ins:romopt-iog-cost-0}),
         $\tau = 10$ (\ref{line:ins:romopt-iog-cost-1}),
         $\tau = 100$ (\ref{line:ins:romopt-iog-cost-2}),
				 $\tau = \infty$ (\ref{line:ins:romopt-iog-cost-infty}).
         The baseline methods used for comparison: fixed, $5$ level
         isotropic sparse grid, SG-ISO, (\ref{line:ins:stochiso}) and
         the dimension-adaptive sparse grid method, SG-TR,
				 (\ref{line:ins:kouri14}). The cost is defined in
				 Eq.~\eqref{eqn:ins:cost}.}
 \label{fig:ins:obj-vs-cost}
\end{figure}
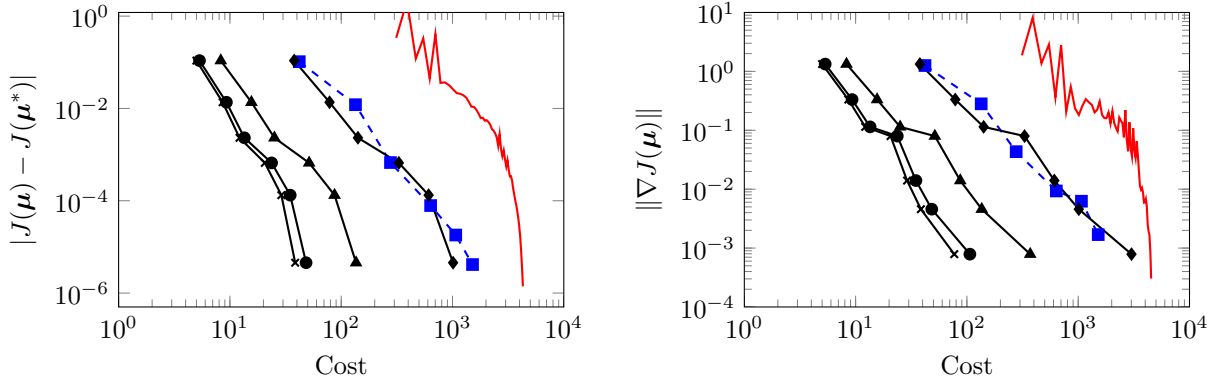
}

\section{Conclusion}
This work proposes and demonstrates the merits of an efficient method for
solving optimization problems constrained by large-scale nonlinear systems of
equations with uncertain parameters
by combining state-of-the-art methods for gradient-based optimization,
stochastic collocation, and the efficient approximation of large-scale systems.
The proposed method approximates the objective function, which comprises an integral
over the stochastic space where the integrand depends
on the solution of the governing system of equations, by combining two
approximation techniques: dimension-adaptive sparse grids and
projection-based reduced-order models. The reduced-order models 
generate inexpensive approximate solutions of the governing system, and sparse grids produce efficient quadrature schemes to integrate the
objective function and its gradient. The proposed approximation model is
used to define the trust-region model problem in the context of a trust-region
method that allows for inexact objective and gradient evaluations, provided
the model objective and gradient are equipped with weak error estimators.
We derive error estimators for our proposed sparse-grid/reduced-basis 
model and introduce dimension-adaptive, greedy algorithms to refine the
sparse grid and reduced basis until the trust-region convergence requirements
are satisfied. Minimum-residual primal and adjoint reduced-order models ensure
these algorithms will terminate and return a sparse grid and reduced basis
that ensure convergence. Unlike existing methods that combine reduced-order
models and sparse grids for stochastic optimization in the context of
linear, elliptic PDEs \cite{chen2015new}, our approach applies to general
nonlinear systems and general quantities of interest with sufficient
regularity; further, the method is guaranteed to be globally convergent to the
solution of the true optimization problem under mild regularity assumptions.
Numerical experiments on a  model problem from optimal stochastic flow control
demonstrated the ability of the proposed method to significantly outperform
alternative approaches in terms convergence as a function of computational cost.
The generality of the method and its promising performance on the
stochastic flow control problem suggest it may be a useful tool to perform
optimization under uncertainty of large-scale problems that commonly arise in
engineering, science, or medical applications, where stochastic optimization is
currently considered too computationally demanding to be practical. An interesting
research direction is to investigate specific relevant applications where the proposed
method will reduce the computational burden of optimization under uncertainty to the
point where it becomes feasible.
Promising future methodological research directions include the extension of the
proposed method to non-smooth risk measures that take into account the semivariance or extreme events, incorporation of hyperreduction to reduce the complexity
of evaluating the nonlinear equation, and {localization} of the reduced basis in the
stochastic space, which will eliminate the need to construct a single basis that must
be accurate over the entire stochastic space.

\appendix

\section{Residual-based error bounds}
\label{app:resbnd}
In this section, we consider a system of nonlinear equations with
sufficient regularity and derive error bounds on the its solution and
quantity of interest in terms of the residual. The nonlinear system
of equations is defined by the mapping
$\func{\dpderes}{\Rbb^\nU\times\Rbb^{\nz}}{\Rbb^\nU}$ and the quantity
of interest is the mapping
$\func{f}{\Rbb^\nU\times\Rbb^{\nz}}{\Rbb}$. The required regularity
and boundedness assumptions imposed on these mappings are stated in
Assumptions~\ref{assume:res}--\ref{assume:qoi}.

\begin{assume} \label{assume:res}
 Consider open, bounded subsets $\Ucal \subset \Rbb^\nU$ and
 $\Zcal \subset \Rbb^{\nz}$. We assume the following:
 \begin{enumerate}
  \item $\func{\dpderes}{\Rbb^\nU\times\Rbb^{\nz}}{\Rbb^\nU}$
        is continuously differentiable with respect to both
        arguments on the domain $\Ucal\times\Zcal$.
  \item For every $\dcombvar \in \Zcal$, there is a unique solution
		$\dpdesolSol$ satisfying $\dpderes(\,\dpdesolSol\,,\,\dcombvar) =
		\zerobold$,
        and the set of solutions
        \begin{equation}
         \Ucal_\Zcal \coloneqq
         \{\dpdesolSol \in \Rbb^\nU \mid \dpderes(\dpdesolSol,\,\dcombvar) =
					\zerobold,
           \,\forall \dcombvar \in \Zcal\}
        \end{equation}
        is a bounded set.
  \item \label{assume:drdu}
        The Jacobian matrix
        $\displaystyle{\func{\pder{\dpderes}{\dpdesol}}
                            {\Rbb^\nU\times\Rbb^{\nz}}
                            {\Rbb^{\nU\times\nU}}}$
        is Lipschitz continuous with respect to its first argument on the
        domain $\Ucal\times\Zcal$.
  \item The parameter Jacobian matrix
        $\displaystyle{\func{\pder{\dpderes}{\dparamvar}}
                            {\Rbb^\nU\times\Rbb^{\nz}}
                            {\Rbb^{\nU\times\nmu}}}$
        is Lipschitz continuous with respect to its first argument on the
        domain $\Ucal\times\Zcal$.
  \item \label{assume:r:D}
        The matrix function
        $\func{\Dbm}{\Rbb^\nU\times\Rbb^\nU\times\Rbb^{\nz}}
              {\Rbb^{\nU\times\nU}}$
        defined as
\begin{equation}
 \Dbm:(\dpdesol_1,\,\dpdesol_2,\,\dcombvar) \mapsto
 \int_0^1 \pder{\dpderes}{\dpdesol}
          (\dpdesol_2+t(\dpdesol_1-\dpdesol_2),\,\dcombvar)\,dt
\end{equation}
        for $\dpdesol_1,\,\dpdesol_2 \in \Ucal$, 
        is invertible with a bounded inverse.
 \end{enumerate}
\end{assume}

\begin{assume} \label{assume:qoi}
 Consider bounded subsets $\Ucal \subset \Rbb^\nU$ and
 $\Zcal \subset \Rbb^{\nz}$. We assume the following:
 \begin{enumerate}
  \item $\func{f}{\Rbb^\nU\times\Rbb^{\nz}}{\Rbb}$
        is continuously differentiable with respect to both
        arguments on the domain $\Ucal\times\Zcal$.
  \item \label{assume:f:f}
        $f$ is Lipschitz continuous with respect to its first argument
        on the domain $\Ucal\times\Zcal$.
  \item $\displaystyle{\func{\pder{f}{\dpdesol}}
                            {\Rbb^\nU\times\Rbb^{\nz}}
                            {\Rbb^\nU}}$
        is Lipschitz continuous with respect to its first argument on the
        domain $\Ucal\times\Zcal$.
  \item $\displaystyle{\func{\pder{f}{\dparamvar}}
                            {\Rbb^\nU\times\Rbb^{\nz}}
                            {\Rbb^\nmu}}$
        is Lipschitz continuous with respect to its first argument on the
        domain $\Ucal\times\Zcal$.
 \end{enumerate}
\end{assume}

\begin{proposition} \label{prop:prim-bnd}
Under Assumptions~\ref{assume:res}--\ref{assume:qoi}, for any
$\bar\dpdesol \in \Ucal$ and $\dcombvar \in \Zcal$, there exists a constant
$\kappa > 0$ such that
\begin{equation} \label{eqn:resbnd0}
 \norm{\dpdesolSol - \bar\dpdesol} \leq
 \kappa \norm{\dpderes(\bar\dpdesol,\,\dcombvar)},
\end{equation}
where $\Ucal \subset \Rbb^\nU$ and $\Zcal \subset \Rbb^{\nz}$ are bounded subsets
 and $\dpdesolSol$ is the unique solution satisfying
$\dpderes(\dpdesolSol,\,\dcombvar) = \zerobold$. 
	Furthermore, there exists a constant
$\kappa' > 0$ such that
\begin{equation} \label{eqn:qoibnd0}
 |f(\dpdesolSol,\,\dcombvar) - f(\bar\dpdesol,\,\dcombvar)| \leq
 \kappa'\norm{\dpderes(\bar\dpdesol,\,\dcombvar)}.
\end{equation}

\begin{proof}
Consider $\dcombvar \in \Zcal$ and let $\dpdesolSol \in \Ucal_\Zcal$ be the
unique solution satisfying $\dpderes(\,\dpdesolSol\,,\,\dcombvar) = \zerobold$. 
	Then for any
$\bar\dpdesol \in \Ucal \subset \Rbb^\nU$, we have
\begin{equation}
 \dpderes(\dpdesolSol,\,\dcombvar) - \dpderes(\bar\dpdesol,\,\dcombvar) =
 \Dbm(\dpdesolSol,\,\bar\dpdesol,\,\dcombvar) \cdot (\dpdesolSol-\bar\dpdesol)
\end{equation}
from the definition of $\Dbm$. From the definition of $\dpdesolSol$ and
Assumption~\ref{assume:res}.\ref{assume:r:D}, we have
\begin{equation}
 \dpdesolSol-\bar\dpdesol = -\Dbm(\dpdesolSol,\,\bar\dpdesol,\,\dcombvar)^{-1}
                          \dpderes(\bar\dpdesol,\,\dcombvar).
\end{equation}
The result in (\ref{eqn:resbnd0}) follows directly from the above expression
and the boundedness of $\Dbm^{-1}$
(Assumption~\ref{assume:res}.\ref{assume:r:D}).
The result in (\ref{eqn:qoibnd0}) follows directly from Lipschitz continuity
of $f$ (Assumption~\ref{assume:qoi}.\ref{assume:f:f}) and (\ref{eqn:resbnd0}).
\end{proof}
\end{proposition}

\begin{proposition} \label{prop:dual-bnd}
Under Assumptions~\ref{assume:res}--\ref{assume:qoi}, for any
$\bar\dpdesol \in \Ucal$, $\bar\dadjvar \in \Lambda$, and
$\dcombvar \in \Zcal$, there exist constants $\kappa,\,\tau > 0$ such that
\begin{equation} \label{eqn:resbnd1}
 \norm{\dadjvarSol-\bar\dadjvar} \leq
 \kappa\norm{\dpderes(\bar\dpdesol,\,\dcombvar)} +
 \tau\norm{\dadjres(\bar\dadjvar,\,\bar\dpdesol,\,\dcombvar)}
\end{equation}
where $\Ucal \subset \Rbb^\nU$,
      $\Lambda \subset \Rbb^\nU$, and
      $\Zcal \subset \Rbb^{\nz}$ are bounded subsets,
$\dpdesolSol \in \Ucal_\Zcal$ is the unique solution satisfying
$\dpderes(\,\dpdesolSol\,,\,\dcombvar) = \zerobold$, and
$\dadjvarSol\in\Rbb^\nU$ is the unique solution satisfying
$\dadjres(\,\dadjvarSol\,,\,\dpdesolSol,\,\dcombvar) = \zerobold$. Furthermore, there exist
constants $\kappa',\,\tau' > 0$ such that
\begin{equation} \label{eqn:qoibnd1}
 \norm{\dadjgrad(\dadjvarSol,\,\dpdesolSol,\,\dcombvar) -
       \dadjgrad(\bar\dadjvar,\,\bar\dpdesol,\,\dcombvar)} \leq
 \kappa'\norm{\dpderes(\bar\dpdesol,\,\dcombvar)} +
 \tau'\norm{\dadjres(\bar\dadjvar,\,\bar\dpdesol,\,\dcombvar)}.
\end{equation}

\begin{proof}
Consider $\dcombvar \in \Zcal$ and let $\dpdesolSol \in \Ucal_\Zcal$ be the
unique solution satisfying $\dpderes(\,\dpdesolSol\,,\,\dcombvar) = \zerobold$ and
$\dadjvarSol \in \Rbb^\nU$ be the unique solution satisfying
	$\dadjres(\,\dadjvarSol\,,\,\dpdesolSol,\,\dcombvar) = \zerobold$. Then for any
$\bar\dadjvar \in \Lambda$, the definition of
$\dadjres$ in (\ref{eqn:adjresDef}) gives the following relation
\begin{equation}
 \begin{aligned}
  \dadjvarSol-\bar\dadjvar &= -\pder{\dpderes}{\dpdesol}(\dpdesolSol,\,\dcombvar)^{-T}
   \left[-\pder{f}{\dpdesol}(\dpdesolSol,\,\dcombvar)^T +
         \pder{\dpderes}{\dpdesol}(\dpdesolSol,\,\dcombvar)^T\bar\dadjvar\right] \\
   &= -\pder{\dpderes}{\dpdesol}(\dpdesolSol,\,\dcombvar)^{-T}
       \dadjres(\bar\dadjvar,\,\dpdesolSol,\,\dcombvar).
 \end{aligned}
\end{equation}
From the boundedness of the Jacobian inverse
(Assumption~\ref{assume:res}.\ref{assume:drdu}) and the
above relation, the adjoint error can be bounded by the adjoint residual
evaluated at the \emph{exact} primal solution and adjoint approximation,
i.e., there exists a constant $\kappa_0 > 0$ such that
\begin{equation}
 \norm{\dadjvarSol-\bar\dadjvar} \leq
 \kappa_0 \norm{\dadjres(\bar\dadjvar,\,\dpdesolSol,\,\dcombvar)}.
\end{equation}
The adjoint residual at the \emph{exact} primal solution and adjoint
approximation can be bounded as
\begin{equation}
 \begin{aligned}
  \norm{\dadjres(\bar\dadjvar,\,\dpdesolSol,\,\dcombvar)}
     &\leq
     \norm{\dadjres(\bar\dadjvar,\,\bar\dpdesol,\,\dcombvar)} +
     \norm{\dadjres(\bar\dadjvar,\,\dpdesolSol,\,\dcombvar) -
           \dadjres(\bar\dadjvar,\,\bar\dpdesol,\,\dcombvar)} \\
     &\leq
     \norm{\dadjres(\bar\dadjvar,\,\bar\dpdesol,\,\dcombvar)} +
     \norm{\pder{\dpderes}{\dpdesol}(\dpdesolSol,\,\dcombvar)-
           \pder{\dpderes}{\dpdesol}(\bar\dpdesol,\,\dcombvar)}
     \norm{\bar\dadjvar} +
     \norm{\pder{f}{\dpdesol}(\dpdesolSol,\,\dcombvar) -
           \pder{f}{\dpdesol}(\bar\dpdesol,\,\dcombvar)}
 \end{aligned}
\end{equation}
where we have used the triangle inequality and the definition of
$\dadjres$ in (\ref{eqn:adjres}). Lipschitz continuity of the quantity of
interest and the Jacobian imply the existence of constants
$\kappa_1,\,\kappa_2 > 0$ such that
\begin{equation} \label{eqn:resbnd1-int0}
 \norm{\dadjres(\bar\dadjvar,\,\dpdesolSol,\,\dcombvar)} \leq
 \norm{\dadjres(\bar\dadjvar,\,\bar\dpdesol,\,\dcombvar)} +
 \left(\kappa_1 + \kappa_2\norm{\bar\dadjvar}\right)
 \norm{\dpdesolSol-\bar\dpdesol}.
\end{equation}
The result in (\ref{eqn:resbnd1}) follows directly from the relation in
(\ref{eqn:resbnd1-int0}) and boundedness of $\Lambda$.

From the definition of $\dadjgrad$ in (\ref{eqn:adjgrad}) and Lipschitz
continuity of $\displaystyle{\pder{f}{\dparamvar}}$, there exists a constant
$\kappa_3 > 0$
such that
\begin{equation}
 \begin{aligned}
  \norm{\dadjgrad(\dadjvarSol,\,\dpdesolSol,\,\dcombvar) -
        \dadjgrad(\bar\dadjvar,\,\bar\dpdesol,\,\dcombvar)}
  &\leq
  \norm{\pder{f}{\dparamvar}(\dpdesolSol,\,\dcombvar) -
        \pder{f}{\dparamvar}(\bar\dpdesol,\,\dcombvar)} +
  \norm{\dadjvarSol^T\pder{\dpderes}{\dparamvar}(\dpdesolSol,\,\dcombvar) -
        \bar\dadjvar^T\pder{\dpderes}{\dparamvar}(\bar\dpdesol,\,\dcombvar)} \\
  &\leq
  \kappa_3\norm{\dpdesolSol-\bar\dpdesol} +
  \norm{\dadjvar^T\pder{\dpderes}{\dparamvar}(\dpdesolSol,\,\dcombvar) -
        \bar\dadjvar^T\pder{\dpderes}{\dparamvar}(\bar\dpdesol,\,\dcombvar)},
 \end{aligned}
\end{equation}
which reduces to
\begin{equation}
 \norm{\dadjgrad(\dadjvarSol,\,\dpdesolSol,\,\dcombvar) -
       \dadjgrad(\bar\dadjvar,\,\bar\dpdesol,\,\dcombvar)} \leq
 \kappa_3\norm{\dpdesolSol-\bar\dpdesol} +
 \norm{\bar\dadjvar^T\left(\pder{\dpderes}{\dparamvar}(\dpdesolSol,\,\dcombvar) -
                           \pder{\dpderes}{\dparamvar}(\bar\dpdesol,\,\dcombvar)
                 \right)} +
 \norm{\left(\dadjvarSol-\bar\dadjvar\right)^T
       \pder{\dpderes}{\dparamvar}(\dpdesolSol,\,\dcombvar)}
\end{equation}
from a simple application of the triangle inequality. Lipschitz continuity
and boundedness of $\displaystyle{\pder{\dpderes}{\dparamvar}}$ on
$\Ucal \times \Zcal$ implies the existence of constants
$\kappa_4,\,\kappa_5,\,\tau' > 0$ such that
\begin{equation} \label{eqn:qoibnd1-int0}
 \norm{\dadjgrad(\dadjvarSol,\,\dpdesolSol,\,\dcombvar) -
       \dadjgrad(\bar\dadjvar,\,\bar\dpdesol,\,\dcombvar)} \leq
 (\kappa_4 + \kappa_5\norm{\dadjvarSol})\norm{\dpdesolSol-\bar\dpdesol} +
 \tau'\norm{\dadjvarSol-\bar\dadjvar}.
\end{equation}
The result in (\ref{eqn:qoibnd1}) follows directly from the above relation
(\ref{eqn:qoibnd1-int0}), boundedness of $\Lambda$, and the previous
results in (\ref{eqn:resbnd0}), (\ref{eqn:resbnd1}).
\end{proof}
\end{proposition}

\section*{Acknowledgments}
MJZ's research was supported in part by the Department of Energy
Computational Science Graduate Fellowship and the Luis W.\ Alvarez Postdoctoral
Fellowship by the Director, Office of Science, Office of Advanced Scientific
Computing Research, of the U.S. Department of Energy under Contract No.\
DE-AC02-05CH11231 (MZ). KTC's research was sponsored by Sandia's Advanced
Simulation and Computing (ASC) Verification and Validation (V\&V)
Project \#103723. DPK's research was sponsored by DARPA EQUiPS grant
SNL 014150709.

\bibliographystyle{siamplain}
\bibliography{biblio}

\begin{thebibliography}{10}

\bibitem{alexandrov1998trust}
{\sc N.~M. Alexandrov, J.~E. Dennis~Jr, R.~M. Lewis, and V.~Torczon}, {\em A
  trust-region framework for managing the use of approximation models in
  optimization}, Structural Optimization, 15 (1998), pp.~16--23.

\bibitem{arian2000trust}
{\sc E.~Arian, M.~Fahl, and E.~W. Sachs}, {\em Trust-region proper orthogonal
  decomposition for flow control}, tech. report, DTIC Document, 2000.

\bibitem{artzner1999coherent}
{\sc P.~Artzner, F.~Delbaen, J.-M. Eber, and D.~Heath}, {\em Coherent measures
  of risk}, Mathematical Finance, 9 (1999), pp.~203--228.

\bibitem{barrault2004empirical}
{\sc M.~Barrault, Y.~Maday, N.~C. Nguyen, and A.~T. Patera}, {\em An empirical
  interpolation method: application to efficient reduced-basis discretization
  of partial differential equations}, C. R. Math., 339 (2004), pp.~667--672.

\bibitem{carlberg2011gnat}
{\sc K.~Carlberg, C.~Bou-Mosleh, and C.~Farhat}, {\em Efficient non-linear
  model reduction via a least-squares {P}etrov--{G}alerkin projection and
  compressive tensor approximations}, International Journal for Numerical
  Methods in Engineering, 86 (2011), pp.~155--181.

\bibitem{chen2014weighted}
{\sc P.~Chen and A.~Quarteroni}, {\em Weighted reduced basis method for
  stochastic optimal control problems with elliptic {PDE} constraint}, SIAM/ASA
  Journal on Uncertainty Quantification, 2 (2014), pp.~364--396.

\bibitem{chen2015new}
{\sc P.~Chen and A.~Quarteroni}, {\em A new algorithm for high-dimensional
  uncertainty quantification based on dimension-adaptive sparse grid
  approximation and reduced basis methods}, Journal of Computational Physics,
  298 (2015), pp.~176--193.

\bibitem{chen2013multilevel}
{\sc P.~Chen, A.~Quarteroni, and G.~Rozza}, {\em Multilevel and weighted
  reduced basis method for stochastic optimal control problems constrained by
  {S}tokes equations}, Numerische Mathematik,  (2013), pp.~1--36.

\bibitem{chen2013weighted}
{\sc P.~Chen, A.~Quarteroni, and G.~Rozza}, {\em A weighted reduced basis
  method for elliptic partial differential equations with random input data},
  SIAM Journal on Numerical Analysis, 51 (2013), pp.~3163--3185.

\bibitem{chen2015sparse}
{\sc P.~Chen and C.~Schwab}, {\em Sparse-grid, reduced-basis {B}ayesian
  inversion}, Computer Methods in Applied Mechanics and Engineering, 297
  (2015), pp.~84--115.

\bibitem{chen2016sparse}
{\sc P.~Chen and C.~Schwab}, {\em Sparse-grid, reduced-basis {B}ayesian
  inversion: Nonaffine-parametric nonlinear equations}, Journal of
  Computational Physics, 316 (2016), pp.~470--503.

\bibitem{clenshaw1960method}
{\sc C.~W. Clenshaw and A.~R. Curtis}, {\em A method for numerical integration
  on an automatic computer}, Numerische Mathematik, 2 (1960), pp.~197--205.

\bibitem{conn2000trust}
{\sc A.~Conn, N.~Gould, and P.~Toint}, {\em Trust Region Methods}, vol.~1,
  SIAM, 2000.

\bibitem{eftang2013approximation}
{\sc J.~L. Eftang, M.~A. Grepl, A.~T. Patera, and E.~M. R{\o}nquist}, {\em
  Approximation of parametric derivatives by the empirical interpolation
  method}, Foundations of Computational Mathematics, 13 (2013), pp.~763--787.

\bibitem{gerstner2003dimension}
{\sc T.~Gerstner and M.~Griebel}, {\em Dimension--adaptive tensor--product
  quadrature}, Computing, 71 (2003), pp.~65--87.

\bibitem{heinkenschlossCVAR}
{\sc M.~Heinkenschloss, B.~Kramer, T.~Takhtaganov, and K.~Willcox}, {\em
  Conditional-value-at-risk estimation via reduced-order models}, ACDL
  Technical Report TR-2017-04,  (2017).

\bibitem{heinkenschloss2002analysis}
{\sc M.~Heinkenschloss and L.~N. Vicente}, {\em Analysis of inexact
  trust-region {SQP} algorithms}, SIAM Journal on Optimization, 12 (2002),
  pp.~283--302.

\bibitem{kouri2013trust}
{\sc D.~P. Kouri, M.~Heinkenschloss, D.~Ridzal, and B.~G. van
  Bloemen~Waanders}, {\em A trust-region algorithm with adaptive stochastic
  collocation for {PDE} optimization under uncertainty}, SIAM Journal on
  Scientific Computing, 35 (2013), pp.~A1847--A1879.

\bibitem{kouri2014inexact}
{\sc D.~P. Kouri, M.~Heinkenschloss, D.~Ridzal, and B.~G. van
  Bloemen~Waanders}, {\em Inexact objective function evaluations in a
  trust-region algorithm for {PDE}-constrained optimization under uncertainty},
  SIAM Journal on Scientific Computing, 36 (2014), pp.~A3011--A3029.

\bibitem{legresley2006application}
{\sc P.~A. LeGresley}, {\em Application of Proper Orthogonal Decomposition
  (POD) to Design Decomposition Methods}, PhD thesis, Stanford University,
  2006.

\bibitem{lin2017non}
{\sc Z.~Lin, D.~Xiao, F.~Fang, C.~Pain, and I.~M. Navon}, {\em Non-intrusive
  reduced order modelling with least squares fitting on a sparse grid},
  International Journal for Numerical Methods in Fluids, 83 (2017),
  pp.~291--306.

\bibitem{maute2009reduced}
{\sc K.~Maute, G.~Weickum, and M.~Eldred}, {\em A reduced-order stochastic
  finite element approach for design optimization under uncertainty},
  Structural Safety, 31 (2009), pp.~450--459.

\bibitem{negri2015reduced}
{\sc F.~Negri, A.~Manzoni, and G.~Rozza}, {\em Reduced basis approximation of
  parametrized optimal flow control problems for the stokes equations},
  Computers \& Mathematics with Applications, 69 (2015), pp.~319 -- 336,
  \url{https://doi.org/https://doi.org/10.1016/j.camwa.2014.12.010},
  \url{http://www.sciencedirect.com/science/article/pii/S0898122114006075}.

\bibitem{nocedal2006numerical}
{\sc J.~Nocedal and S.~Wright}, {\em Numerical Optimization}, Springer, 2006.

\bibitem{patera2007reduced}
{\sc A.~T. Patera and G.~Rozza}, {\em Reduced basis approximation and a
  posteriori error estimation for parametrized partial differential equations},
  tech. report, (C) MIT, Massachusetts Institute of Technology, 2007.

\bibitem{peherstorfer2013model}
{\sc B.~Peherstorfer}, {\em Model order reduction of parametrized systems with
  sparse grid learning techniques}, PhD thesis, Technische Universit{\"a}t
  M{\"u}nchen, 2013.

\bibitem{peherstorfer2012model}
{\sc B.~Peherstorfer, S.~Zimmer, and H.-J. Bungartz}, {\em Model reduction with
  the reduced basis method and sparse grids}, in Sparse grids and applications,
  Springer, 2012, pp.~223--242.

\bibitem{royset2017risk}
{\sc J.~O. Royset, L.~Bonfiglio, G.~Vernengo, and S.~Brizzolara}, {\em
  Risk-adaptive set-based design and applications to shaping a hydrofoil},
  Journal of Mechanical Design, 139 (2017), pp.~101403--101403--8,
  \url{http://dx.doi.org/10.1115/1.4037623}.

\bibitem{rozza2008reduced}
{\sc G.~Rozza, D.~Huynh, and A.~T. Patera}, {\em Reduced basis approximation
  and a posteriori error estimation for affinely parametrized elliptic coercive
  partial differential equations}, Archives of Computational Methods in
  Engineering, 15 (2008), pp.~229--275.

\bibitem{steihaug1983conjugate}
{\sc T.~Steihaug}, {\em The conjugate gradient method and trust regions in
  large scale optimization}, SIAM Journal on Numerical Analysis, 20 (1983),
  pp.~626--637.

\bibitem{toint1981towards}
{\sc P.~L. Toint}, {\em Towards an efficient sparsity exploiting {N}ewton
  method for minimization}, in Sparse Matrices and Their Uses, I.S. Duff, ed.,
  Academic Press, New York, 1981, pp.~57--87.

\bibitem{torlo2018stabilized}
{\sc D.~Torlo, F.~Ballarin, and G.~Rozza}, {\em Stabilized weighted reduced
  basis methods for parametrized advection dominated problems with random
  inputs}, SIAM/ASA Journal on Uncertainty Quantification, 6 (2018),
  pp.~1475--1502.

\bibitem{ullmann2014pod}
{\sc S.~Ullmann and J.~Lang}, {\em {POD}-{G}alerkin modeling and sparse-grid
  collocation for a natural convection problem with stochastic boundary
  conditions}, in Sparse Grids and Applications - Munich 2012, J.~Garcke and
  D.~Pfl{\"u}ger, eds., Cham, 2014, Springer International Publishing,
  pp.~295--315.

\bibitem{xiao2015non}
{\sc D.~Xiao, F.~Fang, A.~Buchan, C.~Pain, I.~Navon, and A.~Muggeridge}, {\em
  Non-intrusive reduced order modelling of the navier--stokes equations},
  Computer Methods in Applied Mechanics and Engineering, 293 (2015),
  pp.~522--541.

\bibitem{yang2017algorithms}
{\sc H.~Yang and M.~Gunzburger}, {\em Algorithms and analyses for stochastic
  optimization for turbofan noise reduction using parallel reduced-order
  modeling}, Computer Methods in Applied Mechanics and Engineering, 319 (2017),
  pp.~217 -- 239,
  \url{https://doi.org/https://doi.org/10.1016/j.cma.2017.02.030},
  \url{http://www.sciencedirect.com/science/article/pii/S0045782516314827}.

\bibitem{zahr2016phd}
{\sc M.~J. Zahr}, {\em Adaptive model reduction to accelerate optimization
  problems governed by partial differential equations}, PhD thesis, Stanford
  University, August 2016.

\bibitem{zahr2015progressive}
{\sc M.~J. Zahr and C.~Farhat}, {\em Progressive construction of a parametric
  reduced-order model for {PDE}-constrained optimization}, International
  Journal for Numerical Methods in Engineering, 102 (2015), pp.~1111--1135.

\bibitem{ZZou_DPKouri_WAquino_2016a}
{\sc Z.~Zou, D.~P. Kouri, and W.~Aquino}, {\em An adaptive sampling approach
  for solving {PDE}s with uncertain inputs and evaluating risk}, in AIAA
  SciTech Proceedings, 2016.

\bibitem{ZZou_DPKouri_WAquino_2018a}
{\sc Z.~Zou, D.~P. Kouri, and W.~Aquino}, {\em A locally adapted reduced basis
  method for solving risk-averse {PDE}-constrained optimization problems}, in
  AIAA SciTech Proceedings, 2018.

\end{thebibliography}
\end{document}